\documentclass{amsart}
\usepackage{amsmath,amssymb,graphicx}
\usepackage[colorlinks=true,pdftex]{hyperref}
\theoremstyle{remark}

 \usepackage{a4wide}

\begin{document}

\title{Numerical study of the stability of the Peregrine breather}

\author[C.~Klein]{Christian Klein}
 \address[C.~Klein]{Institut de 
Math\'ematiques de Bourgogne,
Universit\'e de Bourgogne,
9 avenue Alain Savary,
BP 47970 - 21078 Dijon Cedex,
France}
\email{christian.klein@u-bourgogne.fr}
\author[M.~Haragus]{Mariana Haragus}
\address[M.~Haragus]{Laboratoire de Math\'ematiques de Besan\c con, 
Universit\'e de Franche-Comt\'e,
16 route de Gray,
25030 Besan\c con CEDEX,
France}
\email{mharagus@univ-fcomte.fr}

\begin{abstract}
The Peregrine breather is widely discussed as a model for rogue waves 
in deep water. We present here a detailed numerical study of 
perturbations of the Peregrine breather as a solution to the 
nonlinear Schr\"odinger (NLS) equations. We first address the modulational 
instability   of the constant modulus solution to NLS. Then we study 
numerically
localized and nonlocalized perturbations of the Peregrine breather 
in the linear and fully nonlinear setting. It is shown that the 
solution is unstable against all considered perturbations.  
\end{abstract}
\maketitle

\section{Introduction}
The importance of rogue waves  needs hardly be stressed, and the attempts in the mathematical 
modelling of these phenomena have been intensified in recent years. 
Since the governing equations for hydrodynamics are similar to the 
ones from nonlinear optics, also optical rogue waves have been 
experimentally and analytically studied recently,  
see \cite{dudley} for a review. In deep water and nonlinear optics 
the focusing nonlinear Schr\"odinger (NLS) equation 
\begin{equation}
    iu_{t}+u_{xx}+2|u|^{2}u=0
    \label{NLS}.
\end{equation}
for a complex $u(x,t)$ describing the modulation of a carrier wave
is an important model equation. This is because 
its solutions show nonlinear focusing and a modulation 
instability, i.e., self-induced amplitude modulation of a continuous 
wave propagating in a nonlinear medium, with subsequent generation of 
localized structures.     For a discussion of NLS in the context of 
rogue waves see for instance \cite{OOS} and references therein. It is 
fair to say that no complete understanding of rogue waves has yet 
been reached. For complementary approaches see for instance 
\cite{DBU,BPSS}. 

Remarkably the NLS equation has solutions called \emph{breathers} that 
show the features of the popular definition of rogue waves as 
`waves that appear from nowhere and disappear without a trace', see 
for instance \cite{AAT}. The 
most popular nonperiodic breather solution in this context is the 
Peregrine breather \cite{Peregrine}
\begin{equation}
    u_{Per} = \left(1-\frac{4(1+4it)}{1+4x^{2}+16t^{2}}\right)e^{2it}
    \label{peregrine},
\end{equation}
which is an exact solution of the focusing NLS equation. It has the 
property to have constant modulus for $|t|\to\infty$, thus appearing 
from nowhere and disappearing after some time, and for 
$|x|\to\infty$.  The maximum of $|u|$ is reached for $x=t=0$ and is 
three times the asymptotic value. Exact NLS solutions describing 
multiple breathers are given in \cite{dubard}, for an alternative 
approach to generate breather solutions see \cite{kalla}. 
It is reported that the Peregrine solution has been recently observed 
experimentally in hydrodynamics \cite{cha1,cha2}, 
in plasma physics \cite{bailung} and in nonlinear optics 
\cite{kibler}.

For a solution to be of practical relevance in the context of rogue 
waves, it must be sufficiently robust to be observable in a realistic 
environment where the initial conditions are never exactly given by 
the solution (\ref{peregrine}), but are at best random perturbations 
thereof. 
This means one has to study perturbations of Peregrine initial data 
and to investigate whether the resulting perturbed solution to the NLS 
equation is in some sense close to (\ref{peregrine}). This is a 
difficult task mathematically even in a linearized setting, i.e., for 
small amplitude perturbations, since the background solution is not 
constant in time. The question of nonlinear stability is even more out of 
reach. Thus with the currently available techniques, one has to 
address this question numerically. For periodic breathers, this has 
been recently done in \cite{CC}. An 
experimental study of the stability of the Peregrine breather 
appearing in hydrodynamical experiments against wind was described in 
\cite{cha3}.

It is the goal of the present paper to present a detailed numerical 
study of the stability of the 
Peregrine solution. We consider first small perturbations, 
e.g., the evolution of initial data of the form 
$u(x,0)=u_{Per}(x,0)+v(x)$ with $||v(x)||_{\infty}\ll 
1$ according to the NLS equation (\ref{NLS}). Of special interest are 
localized perturbations $v(x)$ in the Schwartz space of smooth rapidly 
decreasing functions $\mathbb{S}(\mathbb{R})$. The second kind of 
perturbation to be studied is of the form $u_{0}(x)=\sigma 
u_{Per}(x,t_{0})$ with $\sigma\sim 1$. This means we will study 
both rapidly decreasing and non decreasing perturbations of the 
Peregrine breather. 
 
The paper is organized as follows: in section 2 we collect some 
mathematical preliminaries needed for the following sections. In 
section 3 we study the absolute instability of the asymptotic state 
of the Peregrine solution, i.e., of a solution with constant modulus. 
In section 4, we present a numerical study of the linearized (on the 
Peregrine background) NLS equation for localized and nonlocalized 
perturbations. In section 5, localized and nonlocalized perturbations 
of the Peregrine breather are studied for the full NLS equation to 
explore the nonlinear stability of the breather. We add some 
concluding remarks in section 6.

\section{Preliminaries}
In this section we collect some mathematical information on the NLS 
equation and linearizations on the background of the Peregrine 
breather needed for the ensuing stability study. We also briefly 
review the numerical techniques to be applied in the following. 

\subsection{NLS and linearized equation}
It was shown by Zakharov and Shabat \cite{ZS} that
the NLS equation is completely integrable, i.e., it has an infinite 
number of conserved quantities the most prominent of which being the 
\emph{mass}, the square of the $L^{2}$ norm of $u$ and the 
\emph{energy}. Since we want to discuss in this note solutions not 
decaying to 0 at infinity but with $\lim_{|x|\to\infty}|u|=\kappa>0$ (i.e., 
solutions which are not in $L^{2}$), we 
consider instead the quantity
\begin{equation}
    \mathcal{E}=\frac{1}{2}\int_{R}^{}\left(|u_{x}|^{2}-|u|^{2}(|u|^{2}-\kappa)\right) dx
    \label{E},
\end{equation}
which is a combination of mass and energy. 
The complete integrability of the NLS equation is the reason why many 
explicit solutions as solitons or breathers are known. 

The NLS equation has also some conformal symmetry which implies that 
with $u(x,t)$ 
also
\begin{equation}
    u^{\sigma}(x,t)=\sigma u(\sigma x,\sigma^{2}t),\quad 
    \sigma\in\mathbb{R},
    \label{sigma}
\end{equation}
is a solution of NLS. This mean that asymptotically nonvanishing NLS 
solutions could be always transformed by a suitable choice of 
$\sigma$ in (\ref{sigma}) to have modulus 1 at infinity.

We first study linearizations of the NLS equation (\ref{NLS}) of the form 
$u=u_{Per}+v$ which leads in lowest order of $v$ to the equation
\begin{equation}
    iv_{t}+v_{xx}+4|u_{Per}|^{2}v+2u_{Per}^{2}\bar{v}=0
    \label{lNLS}.
\end{equation}
It is straight forward to check that this equation has the conserved 
quantity (this is a consequence of the conservation of the $L^{2}$ 
norm for solutions to the NLS equation)
\begin{equation}
    M = \int_{\mathbb{R}}^{}(u_{Per}\bar{v}+\bar{u}_{Per}v)dx
    \label{M}.
\end{equation}
Note that there are different ways to linearize the NLS equation 
(\ref{NLS}). An 
alternative form is to write $u=u_{Per}(1+\tilde{v})$ 
which implies 
\begin{equation}
    i\tilde{v}_{t}+\tilde{v}_{xx}+2(\ln u_{Per})_{x}\tilde{v}_{x}+4|u_{Per}|^{2}\Re 
    \tilde{v}=0,
    \label{lNLS2}
\end{equation}
for NLS (\ref{NLS}) in lowest order of $\tilde{v}$. The disadvantage of equation 
(\ref{lNLS2}) from a numerical point of view 
is that it is singular at the zeros of the Peregrine 
solution for $t=0$. Even for finite small $t$, the fact that one is 
close to a singularity is numerically problematic. Therefore we study 
in this note for linearizations numerically only equation (\ref{lNLS}). However, equation 
(\ref{lNLS2}) is convenient for analytical approaches, see the 
following section.

\subsection{Numerical methods}
We use in this paper spectral methods in $x$, i.e., methods which 
show for a smooth solution an exponential decrease of the numerical error with the number 
of collocation points. For the time integration, we use fourth order 
methods. The accuracy of the numerical computation is controlled via 
conserved quantities. Because of unavoidable numerical errors, these 
are in an actual computation not exactly conserved, but can be used as a valid 
indicator of the reached precision: as shown for instance in 
\cite{etna}, they tend to overestimate the numerical accuracy by 1-2 
orders of magnitude. 

For rapidly decreasing $v(x,\cdot)$ in (\ref{lNLS}), we use 
a Fourier spectral method for the spatial dependence and an 
exponential time differencing (ETD) scheme due to Cox and Matthews 
\cite{CM} for the 
time integration.   The latter scheme uses an exact integration of 
the linearized equation in Fourier space. The modulus of the 
Fourier coefficients in the spatial coordinate indicates the spatial 
resolution since it is known that they decrease exponentially for 
analytic functions. In all studied examples in this paper, they decrease 
to machine precision (here $\approx 
10^{-15}$) during the whole computation which indicates that the 
solution is always fully resolved spatially.

For nonlocalized perturbations  of the Peregrine 
solution (\ref{peregrine}) in the full NLS equation (\ref{NLS}), 
Fourier methods are not appropriate. The algebraic decrease of the 
solution for $|x|\to\infty$ makes it difficult to periodically 
continue the solution as an analytic function within the finite 
numerical precision as can be done with Schwartz functions. To obtain 
spectral convergence also in this context, a novel numerical approach 
has been presented in \cite{BK}: the real line is divided into several 
domains each of which is mapped to $[-1,1]$ with a M\"{o}bius 
transformation. Near infinity we use $1/x$ as a local coordinate 
which allows to cover the whole real line with this approach. Note 
that we use in contrast to \cite{BK} four domains here among which 
just one compactified one (with $1/x$ as a coordinate). This is due 
to the symmetric (in $x$) problem studied here in contrast to 
\cite{BK}. On each domain we use polynomial interpolation to obtain 
on each interval a spectral method. Imposing at the domain boundaries 
the conditions that the solution be $C^{1}$ there implies with the 
analytic properties of the NLS equation (which is a second order PDE 
in $x$) that we obtain a spectral method on the whole real line. For 
the time integration, we use again a fourth order method. 

For details the reader is referred to \cite{etna} and \cite{BK} and 
references given therein. Note that the numerical approaches are 
completely independent for the localized perturbations in the linear 
case and the nonlinear case. Therefore we compare the results in 
these cases showing that they are qualitatively similar until the 
nonlinearity dominates in the latter case. This gives further 
evidence for the validity of the presented results. 

\section{Absolute instability of the asymptotic state}

A problem in the stability study of the Peregrine solution   
is the fact that it depends both 
on $x$ and $t$. Therefore, for the analytical study in this section, we focus on the stability of the asymptotic state for $x\to\infty$ when $u_{Per}\to e^{2it}$. Instabilities of asymptotic states typically induce instabilities of the corresponding exact solutions, however we do not attempt to give a rigorous proof of this fact for the Peregrine solution.

The asymptotic state $e^{2it}$ of the Peregrine solution $u_{Per}$ has constant modulus, so that for our analysis it is more convenient to use the linearized equation 
(\ref{lNLS2}). For complex perturbations of the form $\tilde{v}=\alpha+i\beta$ we obtain the system
\begin{equation}
    \alpha_{t}+\beta_{xx}=0,\quad \beta_{t}-\alpha_{xx}-4\alpha=0
    \label{ab},
\end{equation}
which has constant coefficients.
Our purpose is to study the spectral properties of the matrix-operator defined through
\[
\mathcal L (\alpha,\beta) = (-\beta_{xx},\alpha_{xx}+4\alpha),
\]
and acting in $L^2(\mathbb R)\times L^2(\mathbb R)$. This function space corresponds to a choice of localized perturbations, but the instability results below also hold for the larger class of nonlocalized perturbations. We compute the spectrum, which coincides with the essential spectrum, of the operator $\mathcal L$, and then the absolute spectrum. Absolute spectra, which are always located to the left of essential spectra in the complex plane, allow to distinguish between absolute instabilities, when  perturbations grow in time at each point of the domain, and  convective instabilities, when though overall norms grow in time, perturbations locally decay: growing perturbations are convected away towards infinity (e.g., see  \cite{SS} and the references therein).

Since the operator $\mathcal L$ has constant coefficients, these spectra can be determined from the dispersion relation 
\begin{equation}
\mathcal D(\lambda,\nu) = \lambda^2+\nu^2(\nu^2+4),
   \label{disp},
\end{equation}
obtained by looking for solutions of (\ref{ab}) of the form $(\alpha,\beta)(x,t) = e^{\lambda t+\nu x}(a,b)$, for complex numbers $\lambda,\nu$ and $a,b$.
First, the essential spectrum, which coincides with the full spectrum of the operator $\mathcal L$ in this case, is the set
\[
\sigma = \{\lambda\in\mathbb C \; ;\; \mathcal D(\lambda, ik)=0, \mbox{ for some } k\in\mathbb R\}.
\]
Then
\[
\sigma = \{\lambda\in\mathbb C \; ;\; \lambda^2=k^2(4-k^2), \mbox{ for some } k\in\mathbb R\},
\]
implying $\sigma = i\mathbb R \cup [-2,2]$. In particular, the set $\sigma$ contains values in the open right hand complex plane showing that the asymptotic state of the Peregrine solution is (essentially) spectrally stable.

Next, the absolute spectrum is found by looking at the relative location in the complex plane of the solutions $\nu\in\mathbb C$ of the dispersion relation $\mathcal D(\lambda,\nu)=0$ for $\lambda\in\mathbb C$ \cite{SS}. More precisely, for complex numbers $\lambda$ with large real part, the dispersion relation (\ref{disp}) has two roots $\nu_1^+(\lambda)$, $\nu_2^+(\lambda)$ with positive real parts, $\mathrm{Re}\, \nu_{1}^+(\lambda)\geq\mathrm{Re}\, \nu_{2}^+(\lambda)>0$, and two roots $\nu_1^-(\lambda)$, $\nu_2^-(\lambda)$ with negative real parts, $\mathrm{Re}\, \nu_{1}^-(\lambda)\leq\mathrm{Re}\, \nu_{2}^-(\lambda)<0$. Then in this case, the absolute spectrum of the operator $\mathcal L$ is simply defined by
\[
\sigma_{abs} = \{\lambda\in\mathbb C\; ;\; 
\mathrm{Re}\,\nu_2^+(\lambda) = \mathrm{Re}\,\nu_2^-(\lambda)\}.
\]
The essential instability found above is absolute if there exists $\lambda\in\sigma_{abs}$ with $\mathrm{Re}\,\lambda>0$.
Since the dispersion relation is invariant under the change $\nu\mapsto -\nu$, we have that $\nu_j^-(\lambda) = -\nu_j^+(\lambda)$, for $j=1,2$, implying that $\mathrm{Re}\,\nu_2^+(\lambda) = \mathrm{Re}\,\nu_2^-(\lambda)=0$ when $\lambda\in\sigma_{abs}$. As a consequence, the absolute spectrum coincides with the essential spectrum here, $\sigma_{abs} = \sigma = i\mathbb R \cup [-2,2]$, showing that the instability of the asymptotic state of the Peregrine solution is absolute.

\section{Numerical simulations for the linearized equation}
In this section, we study perturbations of the Peregrine breather in 
the linearized regime. This means we study initial value problems for 
the linear equation (\ref{lNLS}) both for rapidly decreasing initial 
data and for data proportional to the Peregrine solution. We find in 
both cases that the solution grows strongly in time which 
confirms the absolute instability discussed in the previous section.
Because of the time dependent background leading to a 
nonautonomous equation (\ref{lNLS}), the growth in time is not 
exponential.  

\subsection{Rapidly decreasing perturbations}

\begin{figure}[htb!]
   \includegraphics[width=0.7\textwidth]{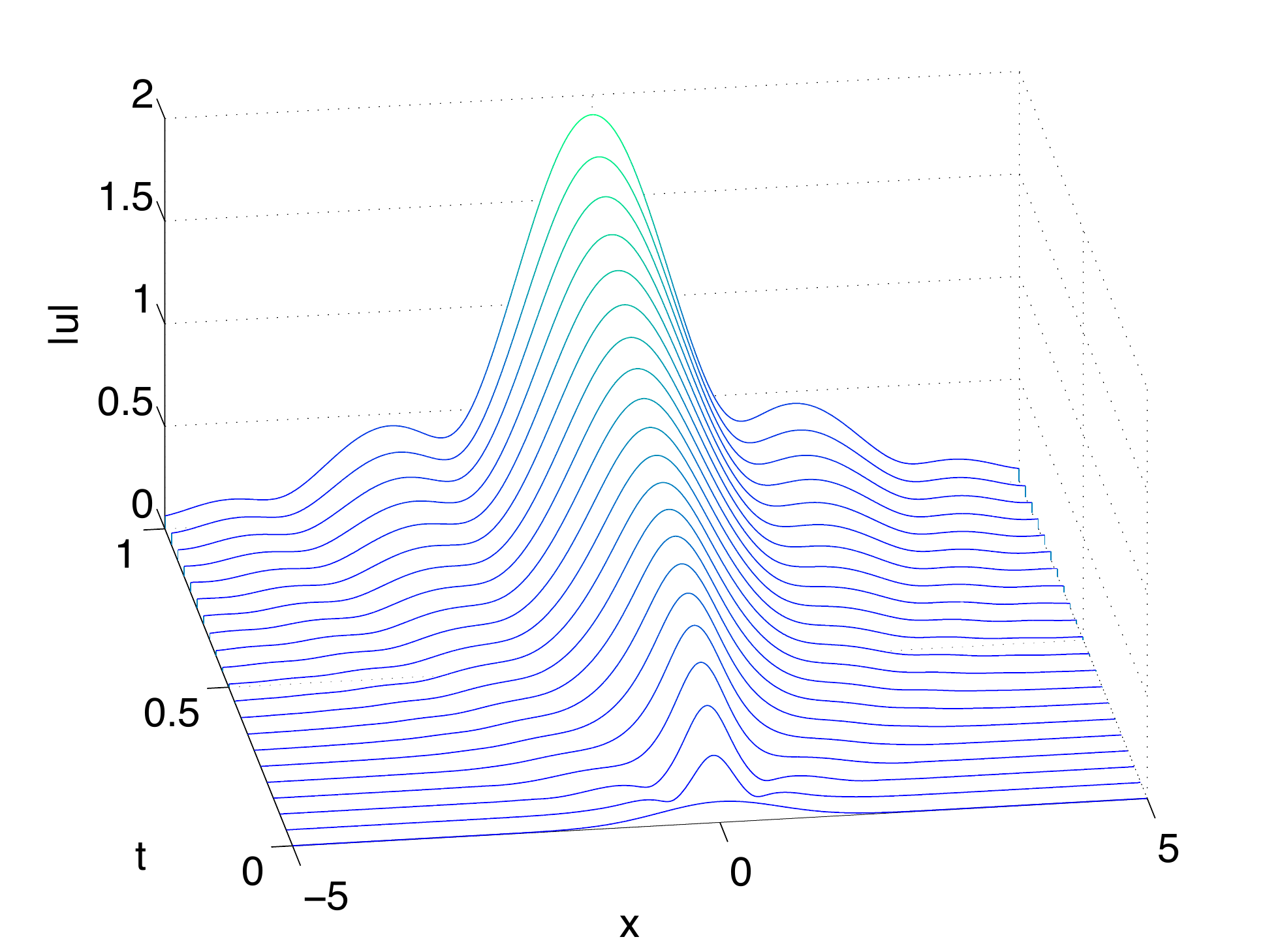}
 \caption{Numerical solution to the linearized NLS equation 
 (\ref{lNLS}) for the initial data $v=0.1e^{-x^{2}}$}.
 \label{nlslinu}
\end{figure}

For the numerical experiments, we concentrate on equation 
(\ref{lNLS}) and the initial data $v=0.1e^{-x^{2}}$. Since these data 
are rapidly decreasing, we use 
as discussed in section 2 a  Fourier spectral method for the spatial dependence.  
The computation is carried out with $N=2^{14}$ 
Fourier modes and $N_{t}=10^{4}$ time steps. The modulus of the 
Fourier coefficients decreases to machine precision 
during the whole computation which indicates that the 
solution is fully resolved spatially. 
The relative change of the conserved quantity 
(\ref{M}) during 
the computation is of the order of $10^{-10}$ which indicates that 
the solution shown in Fig.~\ref{nlslinu} should be accurate at least to the 
order of $10^{-8}$.
The maximum of the solution is visibly quickly growing and becomes 
larger than the modulus of the Peregrine solution itself which 
indicates that the linearized regime is left, in other words that the 
solution is unstable. 
Moreover the initial peak is also being dispersed to 
$\pm\infty$ as can be seen in Fig.~\ref{nlslinu}. 

\subsection{Nonlocalized perturbations}

\begin{figure}[htb!]
   \includegraphics[width=0.7\textwidth]{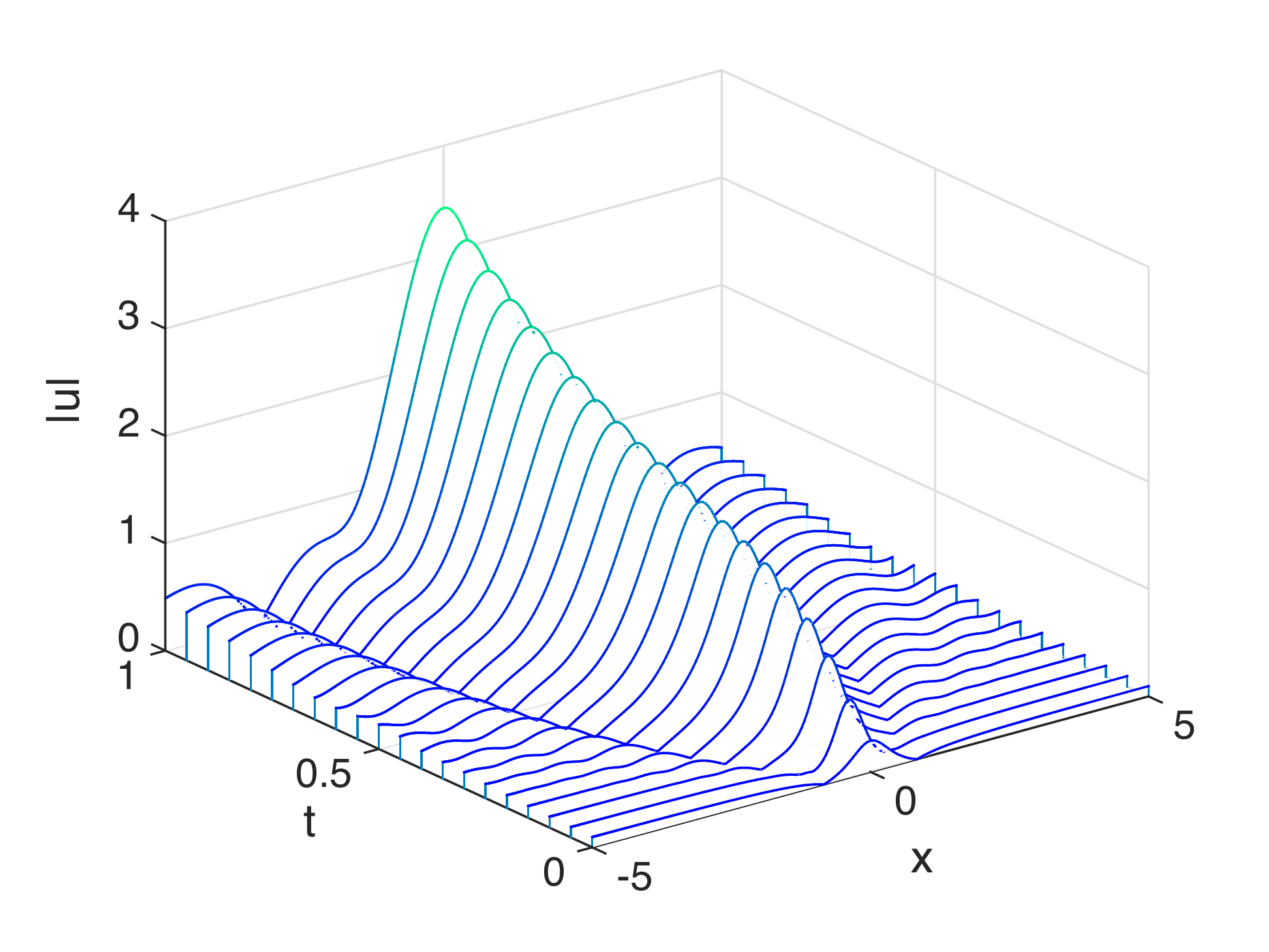}
 \caption{Numerical solution to the linearized NLS equation 
 (\ref{lNLS}) for the initial data $v=0.1u_{Per}$}.
 \label{nlsbreather11lin}
\end{figure}

We consider nonlocalized data  proportional to the breather.
We use the 
same numerical approach as for the fully nonlinear equation and the 
same parameters as detailed there. For the initial data 
$v(x,0)=0.1u_{Per}(x,0)$, we get the solution shown in 
Fig.~\ref{nlsbreather11lin}. As for the localized perturbations in 
Fig.~\ref{nlslinu}, the solution is clearly unstable and quickly 
leaves the linearized regime. 

%

\section{Nonlinear Stability}
In this section, we will study  perturbations of the Peregrine 
breather for the full NLS equation (\ref{NLS}). We again consider 
localized and nonlocalized perturbations to the Peregrine breather as 
in the linearized case and find that the breather is unstable against 
all perturbations. This also applies if the solution is perturbed at an earlier 
time: the resulting NLS solution for these perturbed initial data 
will not be close to the Peregrine breather at $t=0$ or later times.

The numerical approach uses four domains for $x$: domain I with $x\in [-20,-5]$, domain 
II with $x\in[-5,5]$, domain III with $x\in [5,20]$; the compactified 
domain IV is defined via $-1/20<1/x<1/20$. Thus we solve the NLS 
equation on the whole real line. In the respective domains, we use 
$N_{I}=N_{III}=200$, $N_{II}=300$ and $N_{IV}=100$ Chebyshev 
polynomials. We always take the  time step $\Delta_{t}=10^{-3}$. 
The numerical accuracy is controlled via the 
conserved quantity $\mathcal{E}$ in (\ref{E}) in the form of the 
relative 
quantity
\begin{equation}
    \Delta_{E}=1-\frac{\mathcal{E}(t)}{\mathcal{E}(t_{0})},
    \label{DE}
\end{equation}
(because of numerical errors the computed $\mathcal{E}$ will be time 
dependent).

\subsection{Rapidly decreasing perturbations}
\begin{figure}[htb!]
   \includegraphics[width=0.7\textwidth]{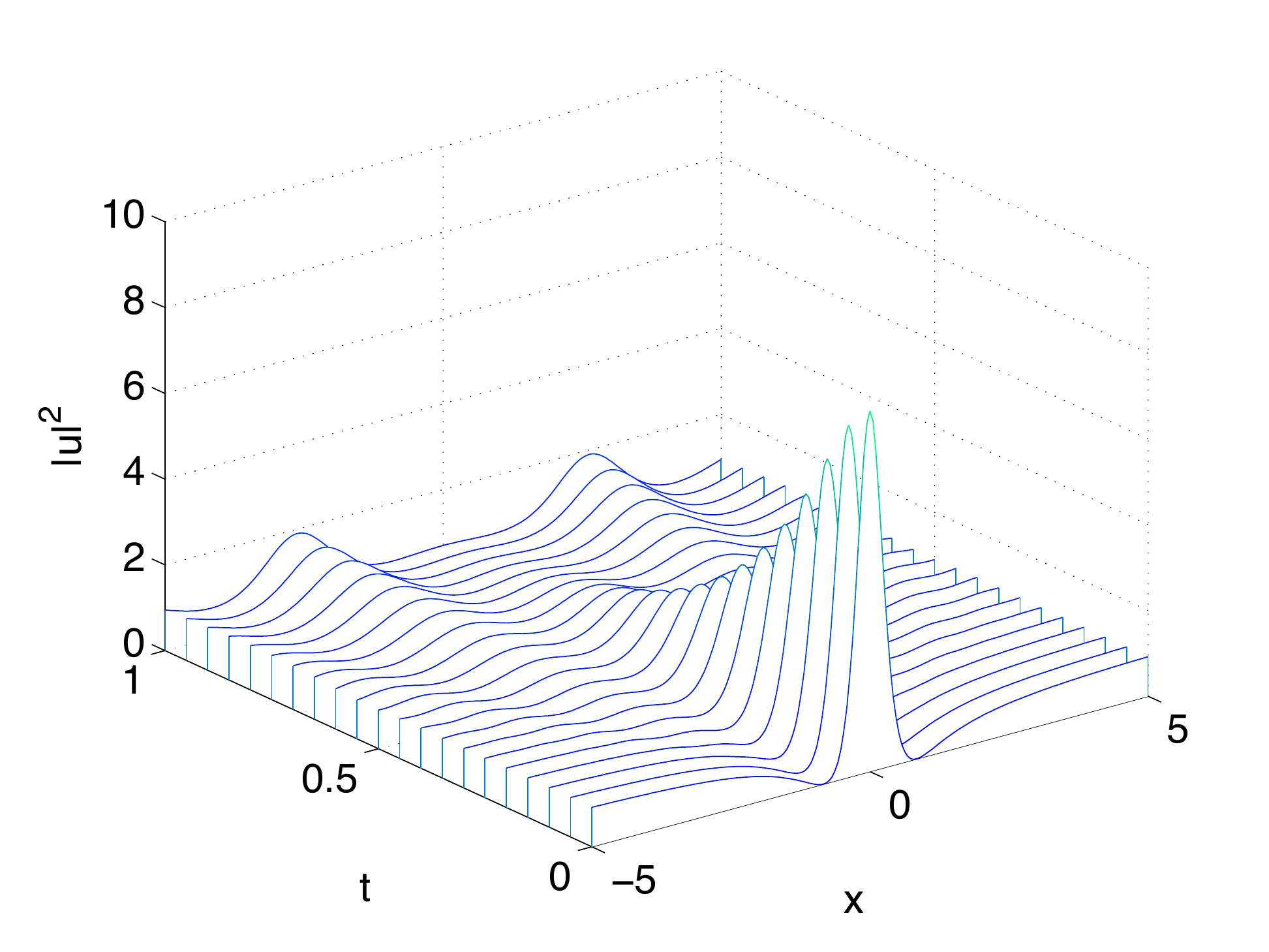}
 \caption{Solution to the NLS equation for the initial data 
 $u_{Per}(x,0)+0.1\exp(-x^{2})$ in dependence of time.}
 \label{breathergauss}
\end{figure}

We first study as in the previous section initial data of the form 
$u(x,0)=u_{Per}(x,0)+0.1\exp(-x^{2})$. The resulting solution can be 
seen in Fig.~\ref{breathergauss} where the initial hump appears to be 
dispersed away to infinity, but in a significantly different way than 
the Peregrine breather. 

\begin{figure}[htb!]
    \includegraphics[width=0.49\textwidth]{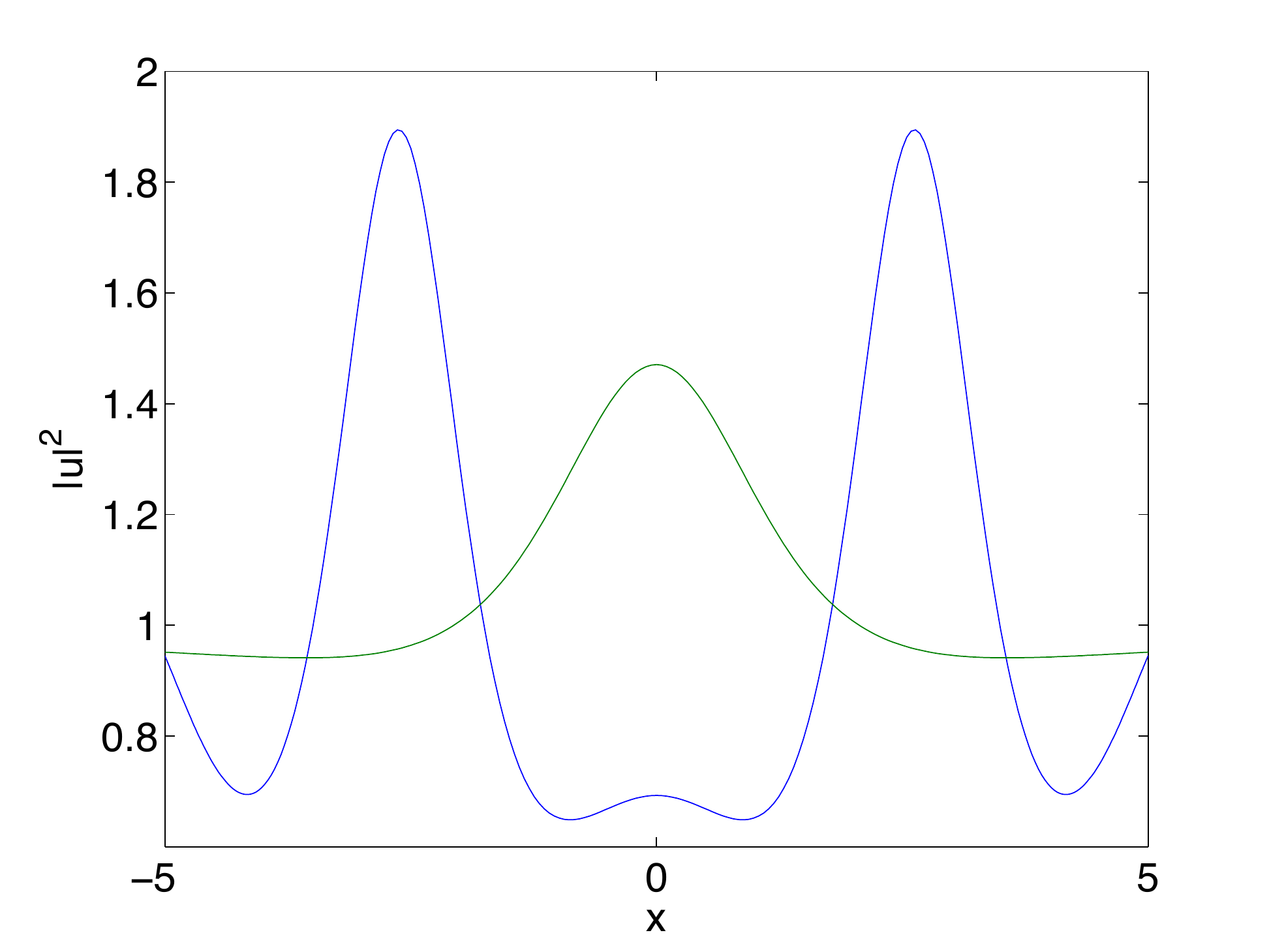}
    \includegraphics[width=0.49\textwidth]{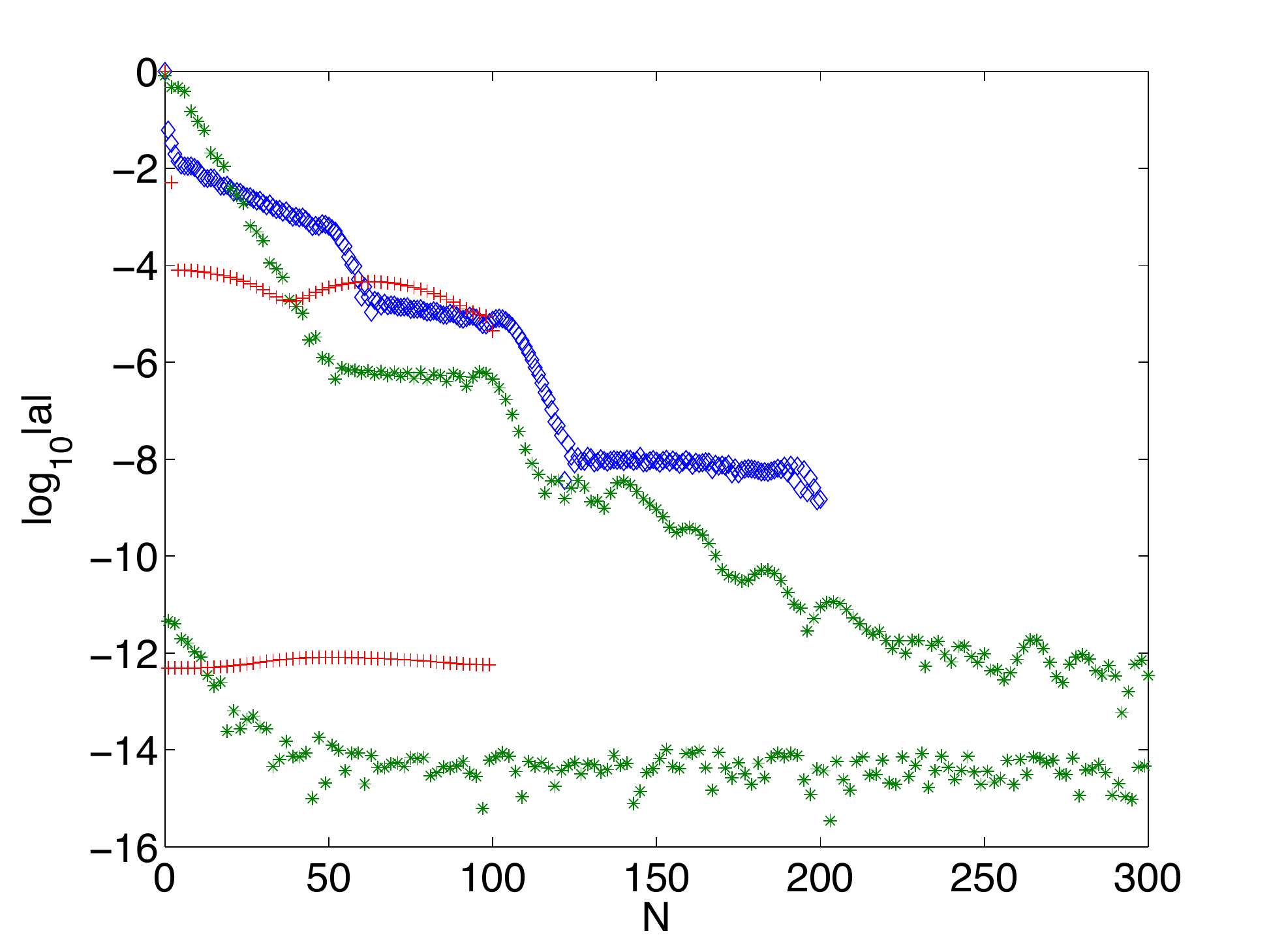}
 \caption{Solution to the NLS equation shown in 
 Fig.~\ref{breathergauss} at $t=1$ in blue and the Peregrine breather 
 (\ref{peregrine}) for $t=1$ in green on the left, and the 
 corresponding Chebyshev coefficients for the numerical 
 solution on the right (in blue for domain I, green for domain II and 
 red for domain III).}
 \label{breathergausscheb}
\end{figure}

This is even more obvious in the left figure of 
Fig.~\ref{breathergausscheb} where the solution at the last computed time 
is shown together with the Peregrine solution (\ref{peregrine}) at the same 
time. Clearly the Peregrine breather is unstable against 
this type of perturbations, and the perturbed solution does not 
stay close to the exact solution. Note that the solution is even at 
the last computed time well resolved spatially as can be seen in the 
right figure of Fig.~\ref{breathergausscheb}, where the Chebyshev 
coefficients at the last computed time are shown. Since 
the initial data are symmetric with respect to the transformation $x\to-x$ and since the 
Schr\"odinger equation preserves parity, the Chebyshev coefficients 
in zone I and III are identical (therefore the coefficients in zone 
III are not shown), and half of the coefficients in zone 
II and IV vanish with numerical precision. The coefficients in zone IV 
still decrease to $10^{-4}$. As discussed in \cite{BK}, 
the slow decrease of the 
Chebyshev coefficients is due to the oscillations with a slowly 
decreasing amplitude for $|x|\to\infty$.  The conserved quantity 
(\ref{DE}) is
$\Delta_{E}\sim 8.8*10^{-3}$. This implies together with the 
resolution in coefficient space shown in 
Fig.~\ref{breathergausscheb} that the solution is computed to at 
least plotting accuracy. 

\begin{figure}[htb!]
    \includegraphics[width=0.7\textwidth]{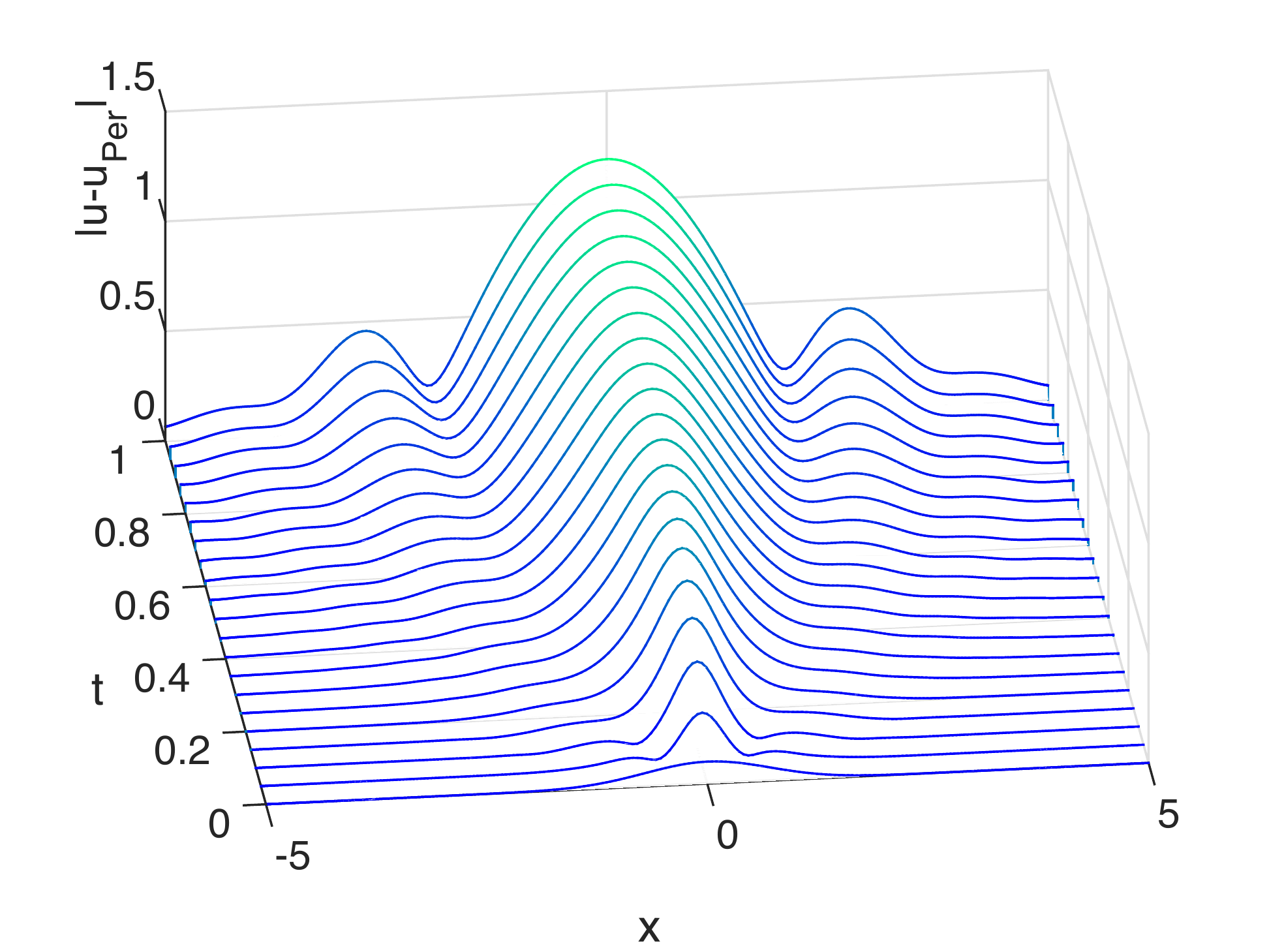}
 \caption{Difference between the solution to the NLS equation shown in 
 Fig.~\ref{breathergauss}  and the Peregrine breather 
 (\ref{peregrine}).}
 \label{nlsbreathergaussdiff}
\end{figure}

In Fig.~\ref{nlsbreathergaussdiff} we show the difference between the 
numerical solution $u$ in Fig.~\ref{breathergauss} and the Peregrine 
solution. This corresponds to the function $v$ in the linearized case 
(\ref{lNLS}), and Fig.~\ref{nlsbreathergaussdiff} thus has to be 
compared to Fig.~\ref{nlslinu}. It is interesting to see that both 
figures are qualitatively similar, but that the nonlinearity has the 
effect to delimit the growing of $|u-u_{Per}|$ for later times and to 
disperse the solution more efficiently. 

\begin{figure}[htb!]
   \includegraphics[width=0.7\textwidth]{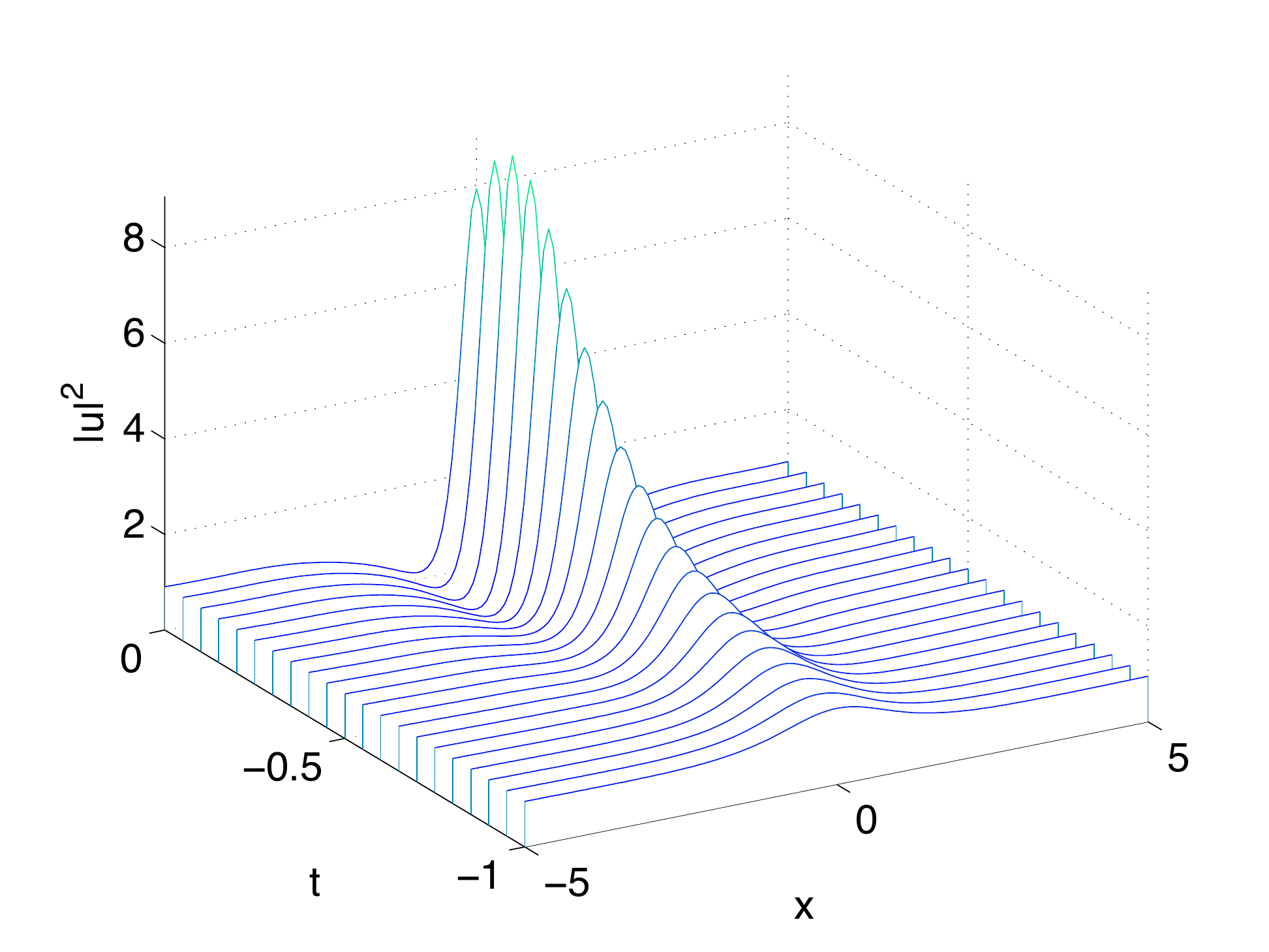}
 \caption{Solution to the NLS equation for the initial data 
 $u(x,-1)=u_{Per}(x,-1)+0.1\exp(-x^{2})$ in dependence of time.}
 \label{breathergausstm1}
\end{figure}

It is an interesting question whether a perturbation of the breather 
at an earlier time still allows to observe a Peregrine type behavior 
at $t\sim 0$ at least in a qualitative sense. To answer this 
question, we consider as initial data 
$u(x,-1)=u_{Per}(x,-1)+0.1\exp(-x^{2})$  in 
Fig.~\ref{breathergausstm1}. It can be recognized  that the solution 
grows to a similar maximal value as the Peregrine solution, but 
reaches it at an earlier time. The focusing effect of the focusing 
NLS equation is clearly visible, but it does not force the solution 
towards the Peregrine breather. 

\begin{figure}[htb!]
    \includegraphics[width=0.49\textwidth]{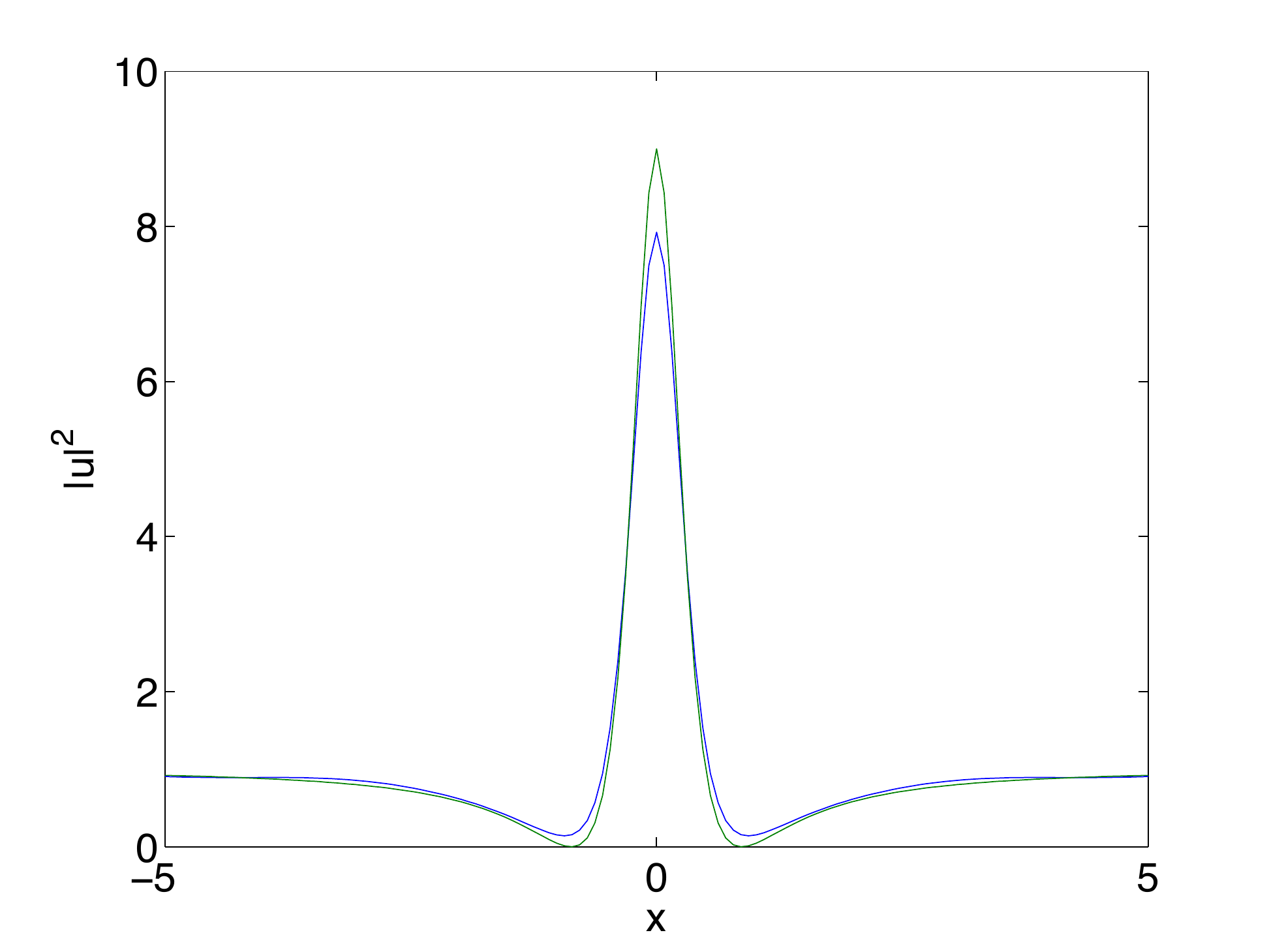}
    \includegraphics[width=0.49\textwidth]{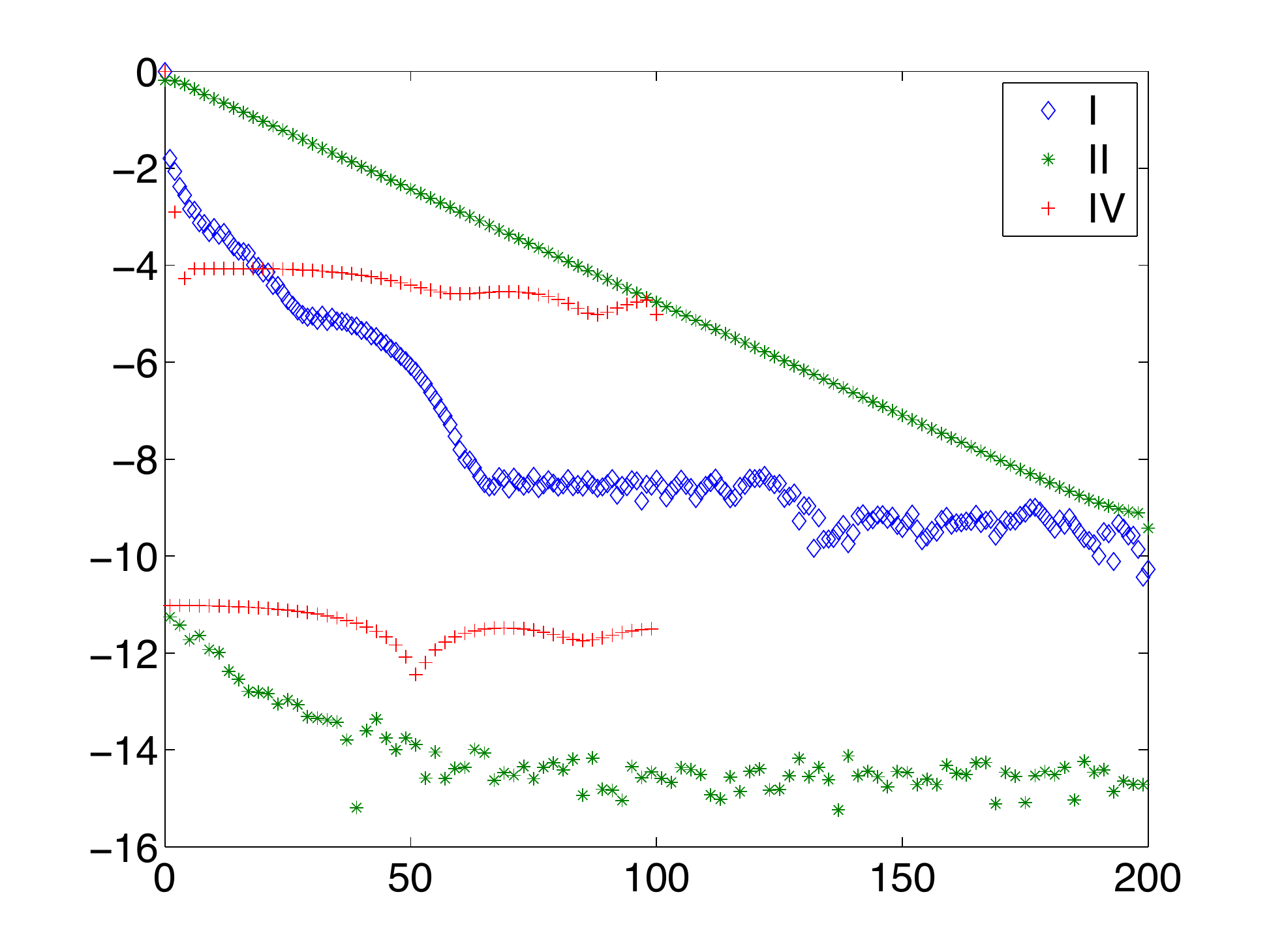}
 \caption{Solution to the NLS equation shown in 
 Fig.~\ref{breathergausstm1} at $t=0$ in blue and the Peregrine breather 
 (\ref{peregrine}) for $t=0$ in green on the left, and the 
 corresponding Chebyshev coefficients for the numerical 
 solution on the right (in blue for domain I, green for domain II and 
 red for domain III).}
 \label{breathergausschebtm1}
\end{figure}

In Fig.~\ref{breathergausschebtm1}, it can be seen that the solution 
at the last recorded time is clearly different from the Peregrine 
solution \ref{peregrine} which again confirms that the solution is 
unstable. The Chebyshev coefficients at the last recorded time can be 
seen in Fig.~\ref{breathergausschebtm1} on the right showing that the 
solution is well resolved. 

\subsection{Nonlocalized perturbations}
In this subsection we study nonlocalized perturbations of the 
breather of the form $u(x,t_{0})=\sigma u_{Per}(x,t_{0})$, where 
$\sigma\sim 1$. The parameters for the  numerical computation 
are as in the previous subsection. 

\begin{figure}[htb!]
   \includegraphics[width=0.7\textwidth]{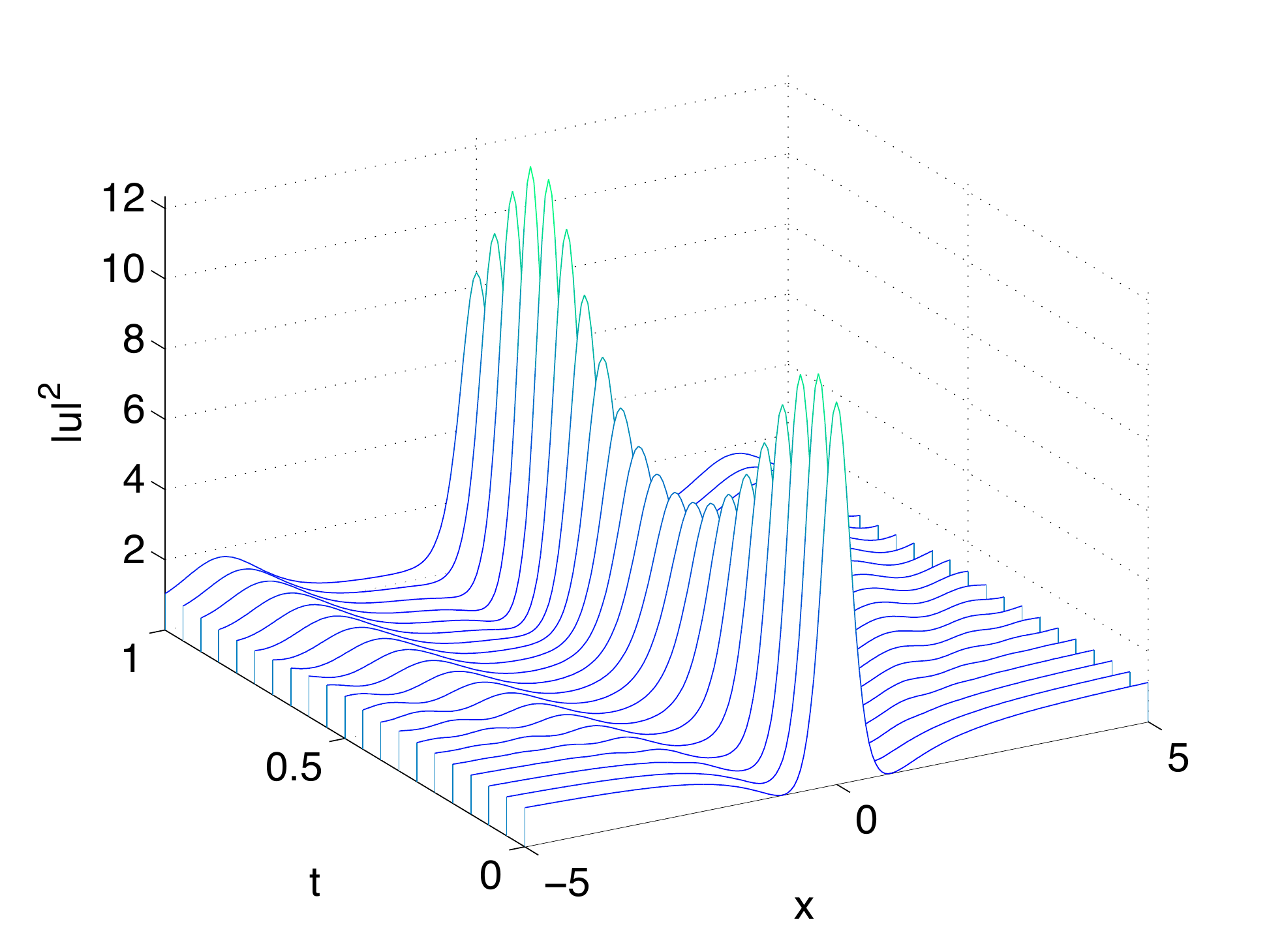}
 \caption{Solution to the NLS equation for the initial data 
 $u(x,0)=1.1u_{Per}(x,0)$ in dependence of time.}
 \label{breather11}
\end{figure}

In Fig.~\ref{breather11} we show the solution for the initial data 
$u(x,0) = 1.1u_{Per}(x,0)$ in dependence of time. It can be seen that 
the solution keeps growing despite being initially already above the 
global maximum of the Peregrine solution before decreasing, but in 
contrast to the latter not monotonically. Instead it grows again 
before being eventually dispersed.
Note that the quantity $\mathcal{E}$ in (\ref{E}) is positive in 
this case.

\begin{figure}[htb!]
    \includegraphics[width=0.49\textwidth]{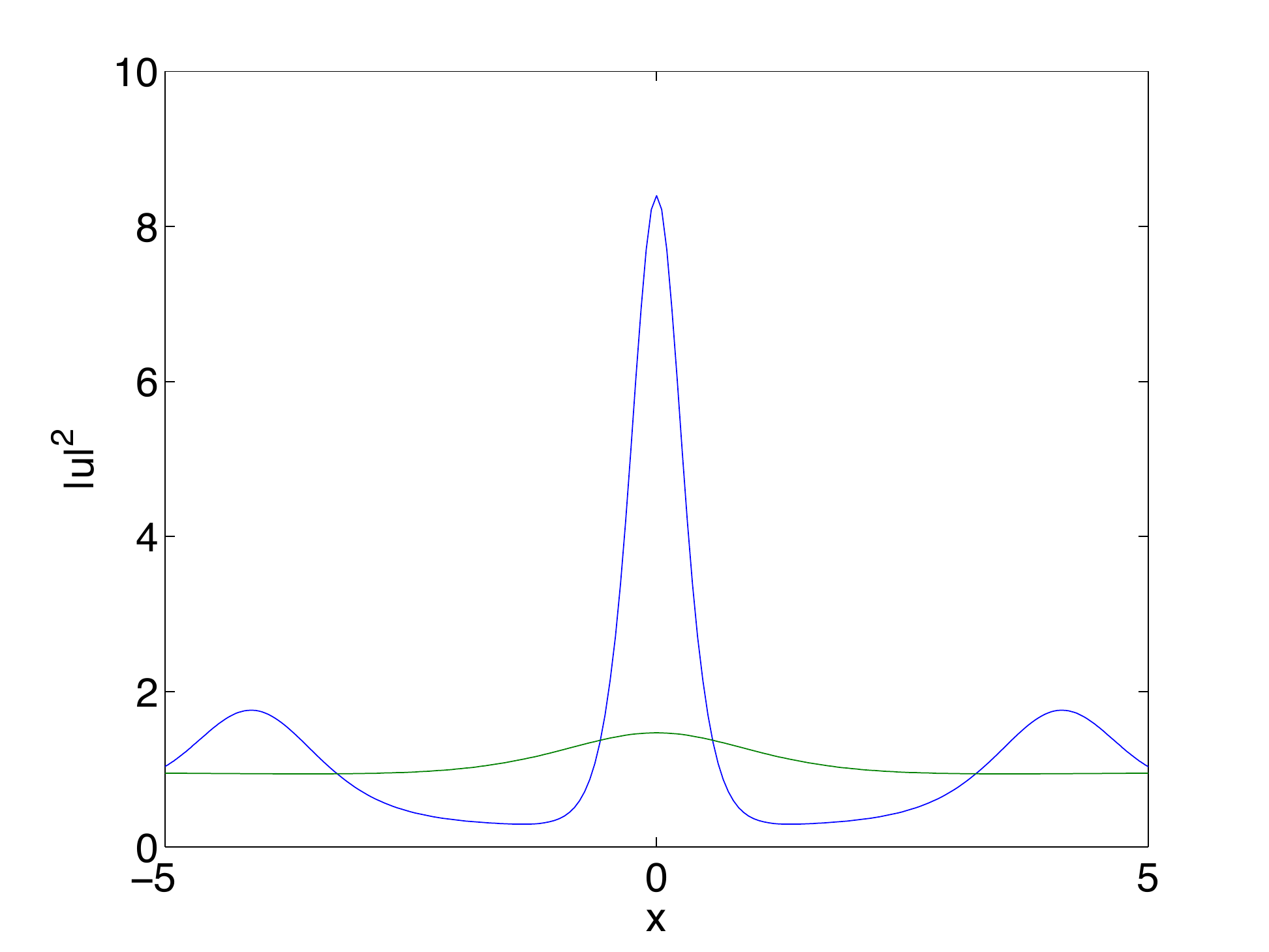}
    \includegraphics[width=0.49\textwidth]{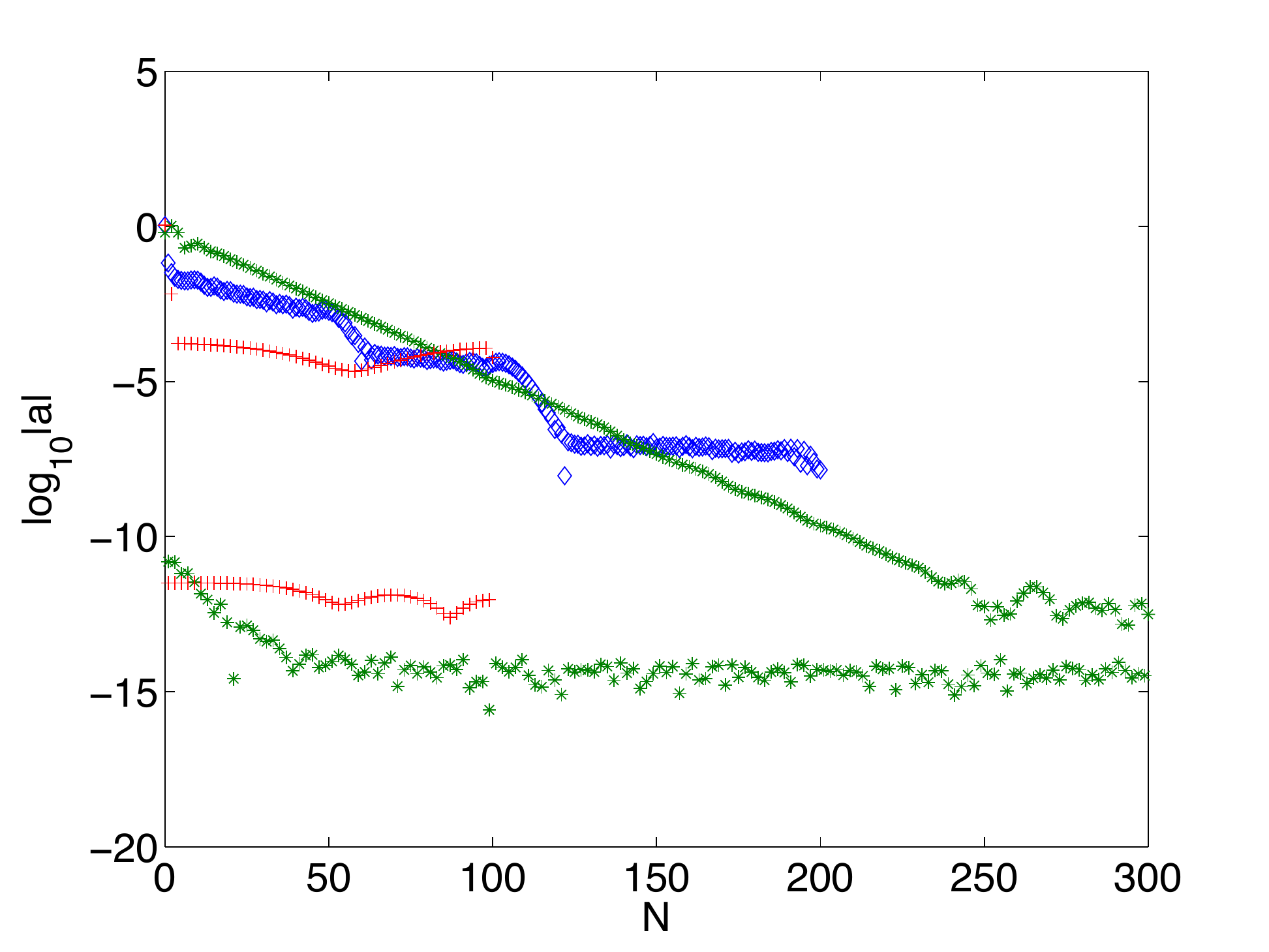}
 \caption{Solution to the NLS equation shown in 
 Fig.~\ref{breather11} at $t=1$ in blue and the Peregrine breather 
 (\ref{peregrine}) for $t=1$ in green on the left, and the 
 corresponding Chebyshev coefficients for the numerical 
 solution on the right (in blue for domain I, green for domain II and 
 red for domain III).}
 \label{breather11cheb}
\end{figure}

It is obvious, see also the left figure in Fig.~\ref{breather11cheb},
that the solution at the final recorded time is not 
close to the Peregrine solution at the same time. Clearly the 
Peregrine solution is also nonlinearly unstable against unbounded 
perturbations proportional to the solution itself. The right figure 
in Fig.~\ref{breather11cheb} shows that the solution is again well 
resolved spatially even at the final time of computation. For the quantity (\ref{DE}) we have 
$\Delta_{E}<8.15*10^{-4}$ during the whole computation. 

\begin{figure}[htb!]
    \includegraphics[width=0.7\textwidth]{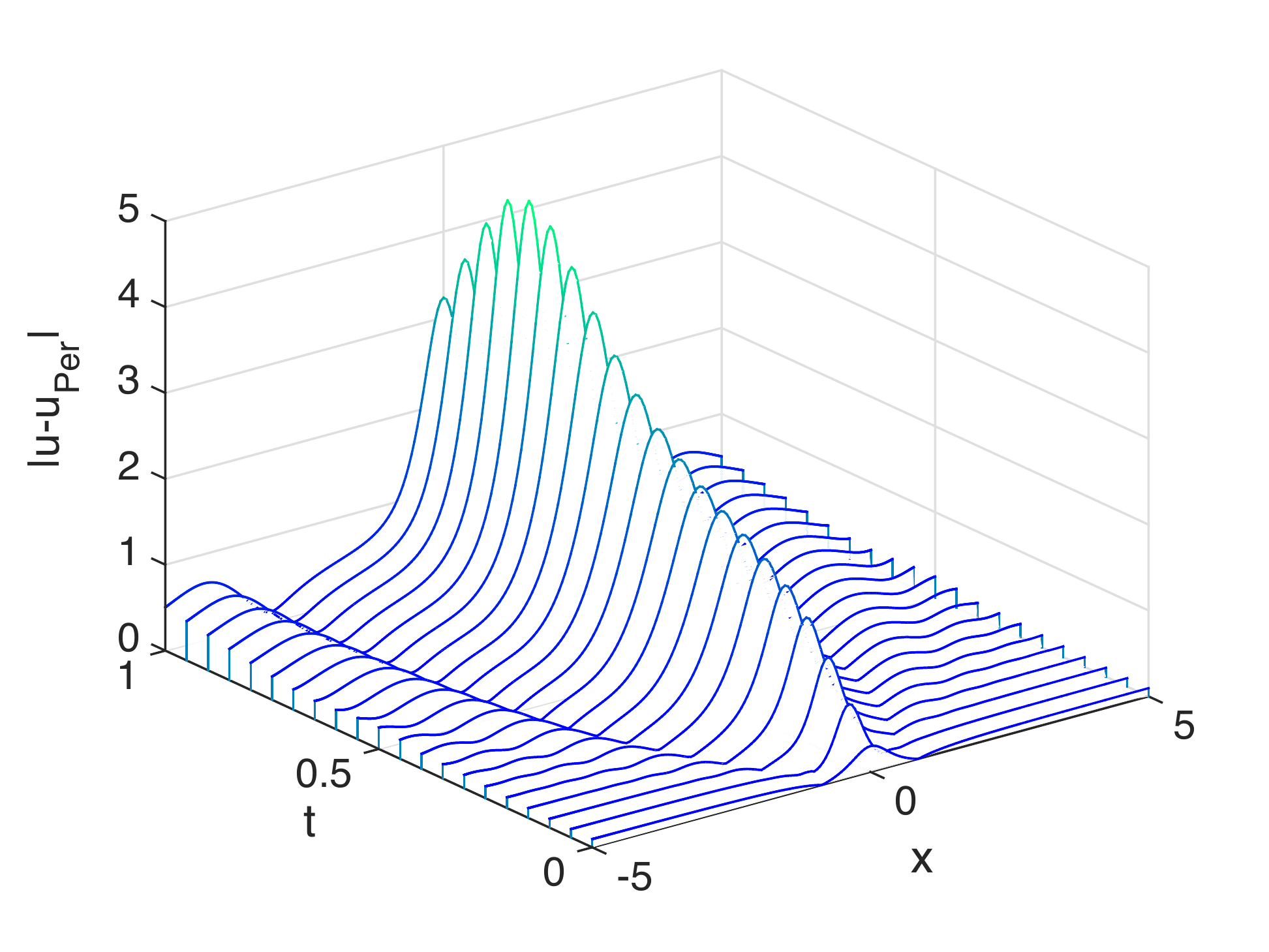}
 \caption{Difference of the solution to the NLS equation shown in 
 Fig.~\ref{breather11} and the Peregrine breather 
 (\ref{peregrine}).}
 \label{nlsbreather11diff}
\end{figure}

In Fig.~\ref{nlsbreather11diff} we show the difference of the 
numerical solution in Fig.~\ref{breather11} and the Peregrine 
breather in dependence of $t$. This can be compared to the solution 
of the linearized equation (\ref{lNLS}) in 
Fig.~\ref{nlsbreather11lin}. It can be seen that the solution to the 
linearized equation (\ref{lNLS}) grows continually, whereas the 
solution to the full NLS equation reaches a maximum and decreases 
then. Note however the qualitative similarity at earlier times. 

\begin{figure}[htb!]
   \includegraphics[width=0.7\textwidth]{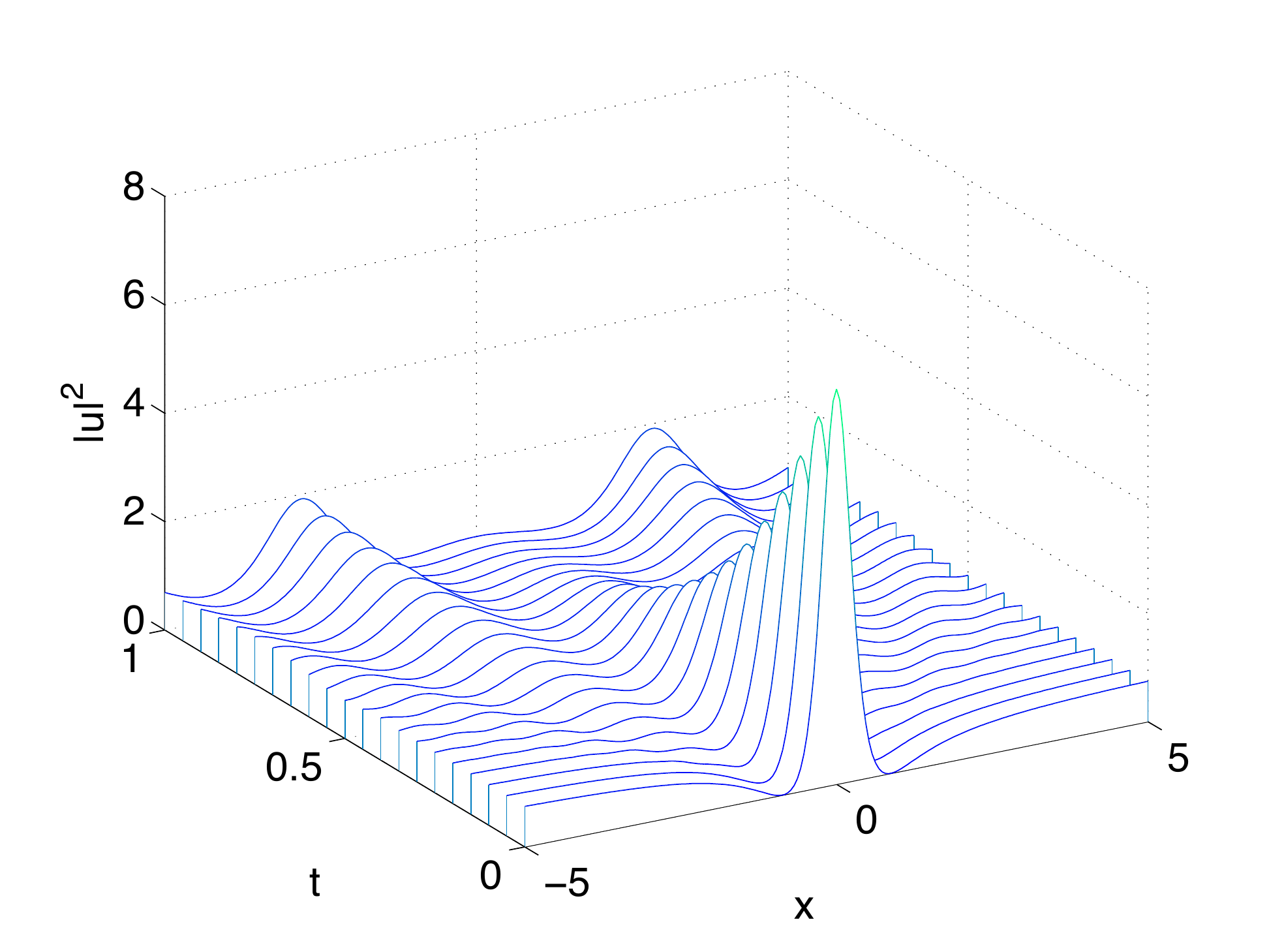}
 \caption{Solution to the NLS equation for the initial data 
 $u(x,0)=0.9u_{Per}(x,0)$ in dependence of time.}
 \label{breather09}
\end{figure}

A similar perturbation as in Fig.~\ref{breather11} but of positive 
$\mathcal{E}$ is considered in Fig.~\ref{breather09} where the initial data 
are $u(x,0)= 0.9u_{Per}(x,0)$. Here the solution never reaches the 
maximum of the Peregrine solution. Instead it decreases in modulus 
with time, but also not as the latter solution, but by radiating the 
mass towards infinity in the form of dispersive oscillations. 

\begin{figure}[htb!]
    \includegraphics[width=0.49\textwidth]{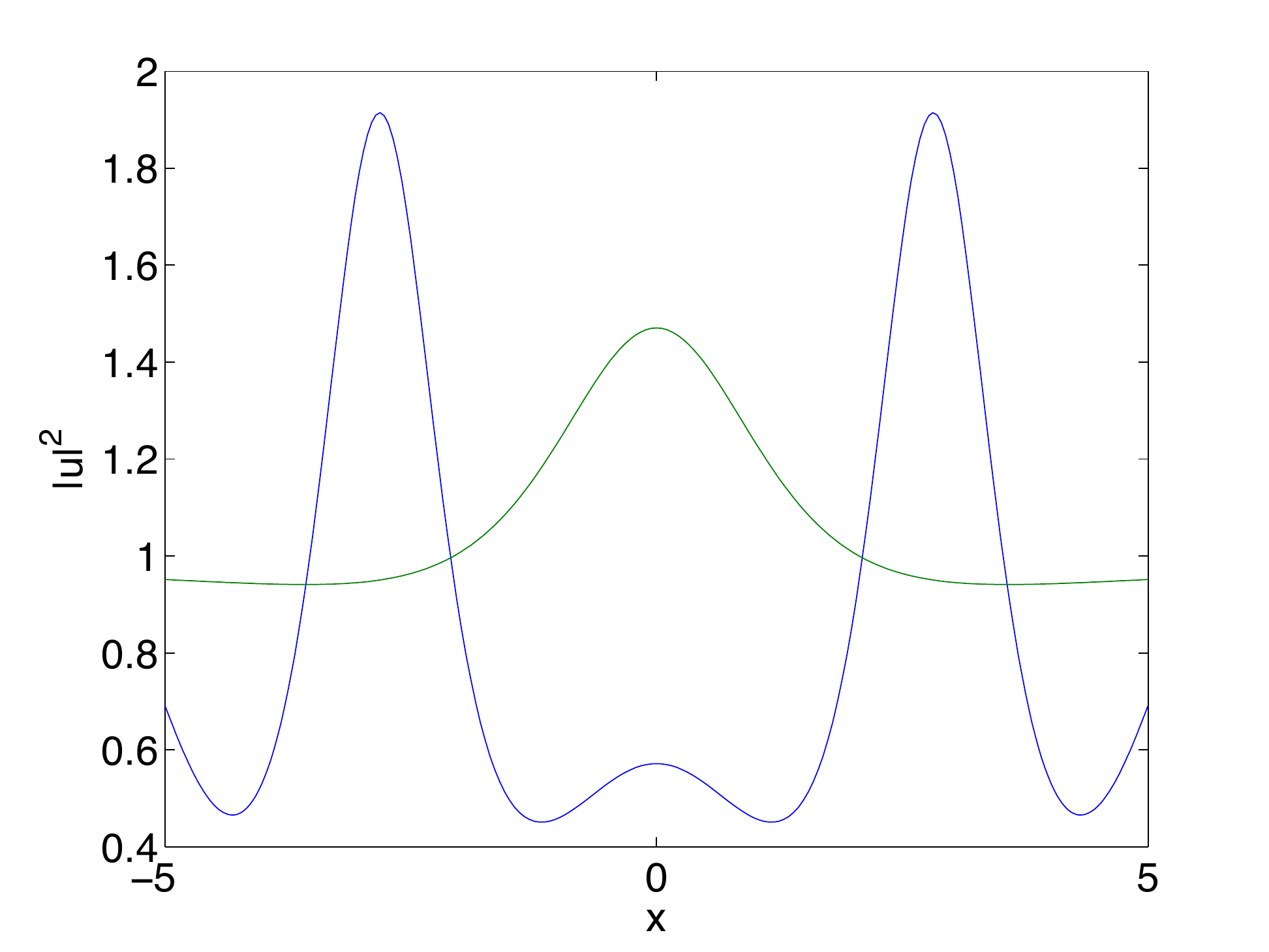}
    \includegraphics[width=0.49\textwidth]{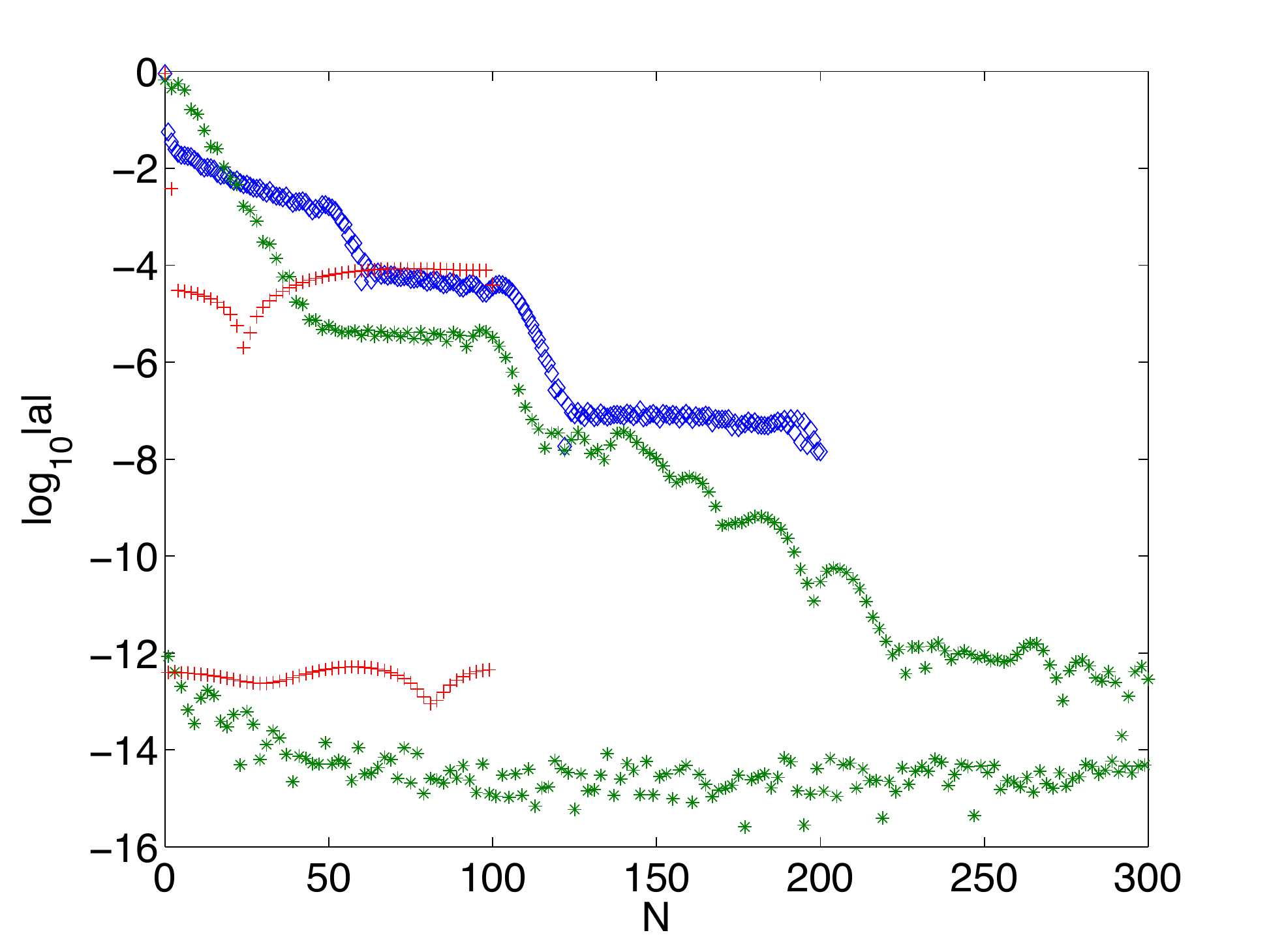}
 \caption{Solution to the NLS equation shown in 
 Fig.~\ref{breather09} at $t=1$ in blue and the Peregrine breather 
 (\ref{peregrine}) for $t=1$ in green on the left, and the 
 corresponding Chebyshev coefficients for the numerical 
 solution on the right (in blue for domain I, green for domain II and 
 red for domain III).}
 \label{breather09cheb}
\end{figure}

Fig.~\ref{breather09cheb} shows that the solution is again in no 
sense close to the Peregrine solution (\ref{peregrine}). The 
Chebyshev coefficients in the same figure decrease again to the order 
of $10^{-4}$. During the whole computation we have 
$\Delta_{E}<4.43*10^{-5}$.

\begin{figure}[htb!]
   \includegraphics[width=0.7\textwidth]{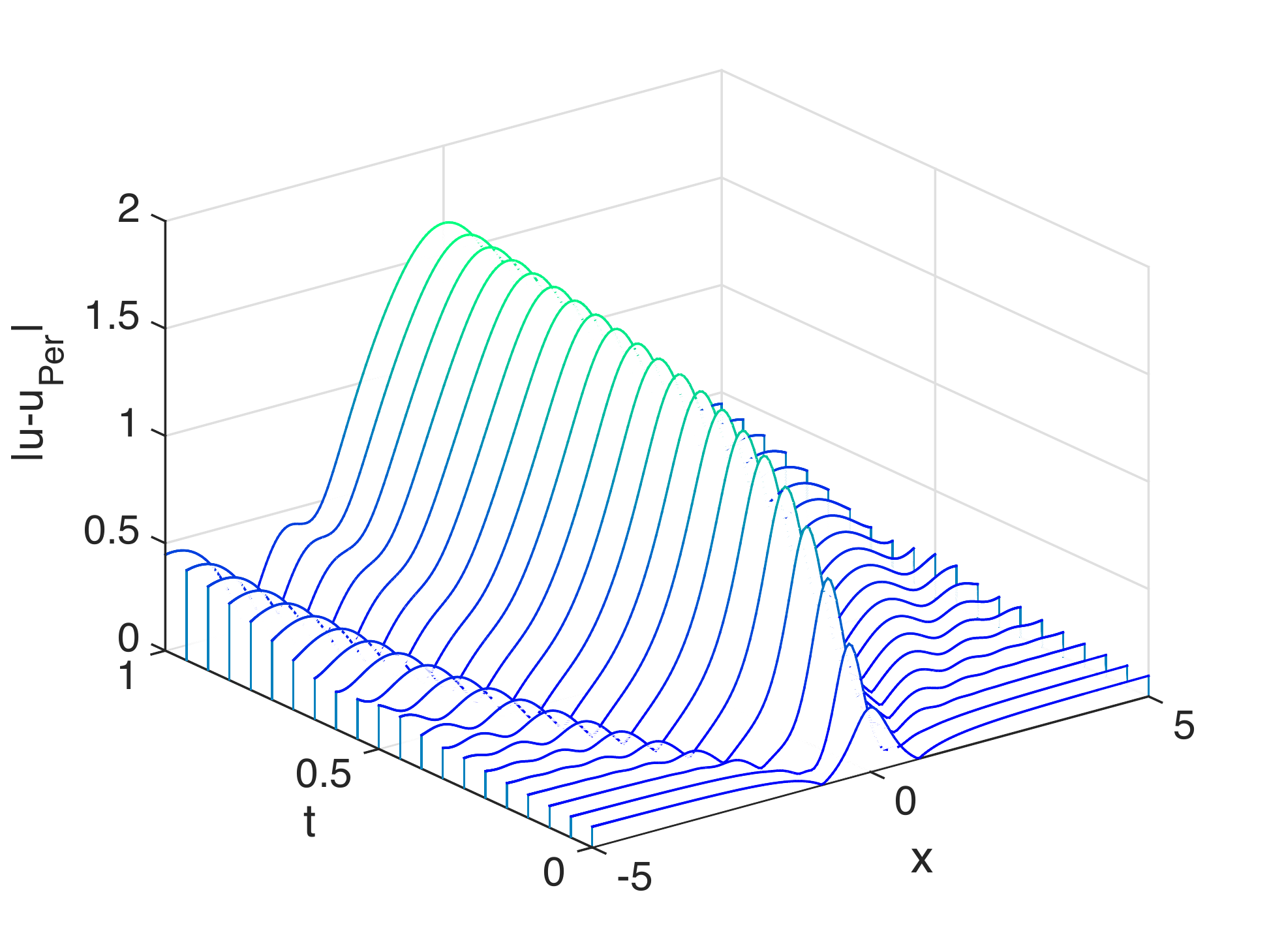}
 \caption{Difference of the solution to the NLS equation in 
 Fig.~\ref{breather09} and the Peregrine breather.}
 \label{nlsbreather09diff}
\end{figure}

In Fig.~\ref{nlsbreather09diff} we show the difference of the 
solution in Fig.~\ref{breather09} and the Peregrine breather which 
can be compared to the solution of the linearized equation 
(\ref{lNLS}) in Fig.~\ref{nlsbreather11lin}. Note that the solution 
of the linearized equation (\ref{lNLS}) is identical for the cases 
corresponding to Fig.~\ref{breather11} and Fig.~\ref{breather09}. 
Clearly this is not the case for the solutions to the full NLS 
equation, but there is good agreement for early times. For larger 
times, the difference in Fig.~\ref{nlsbreather09diff} grows 
considerably more slowly than the solution in 
Fig.~\ref{nlsbreather11lin}. 

\begin{figure}[htb!]
   \includegraphics[width=0.7\textwidth]{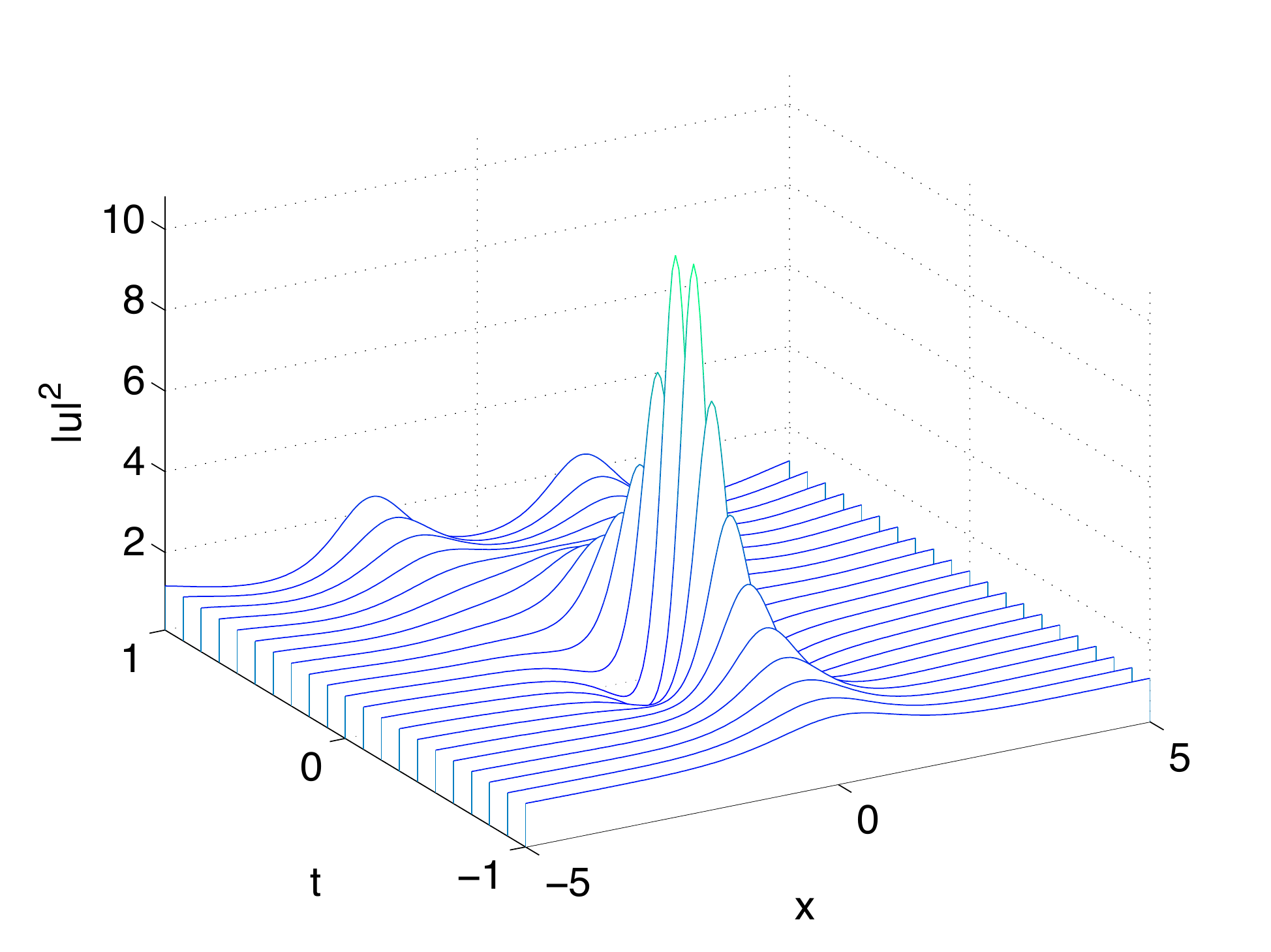}
 \caption{Solution to the NLS equation for the initial data 
 $u(x,-1)=1.1u_{Per}(x,-1)$ in dependence of time.}
 \label{breather11tm1}
\end{figure}

As in the previous subsection, we also address the question whether 
perturbations of the above kind at an earlier time still allow to 
qualitatively observe characteristic features of the Peregrine 
breather (\ref{peregrine}) at later times. In Fig.~\ref{breather11tm1} we consider 
the solution for initial data of the form 
$u(x,-1)=1.1u_{Per}(x,-1)$. 
In this case $\mathcal{E}$ is again positive, and a higher 
maximum is reached at an earlier time than for the Peregrine breather, 
though the solution is qualitatively of a similar form.  However, the 
solution is dispersed away to infinity in a completely different 
manner than for the Peregrine breather. 

\begin{figure}[htb!]
    \includegraphics[width=0.49\textwidth]{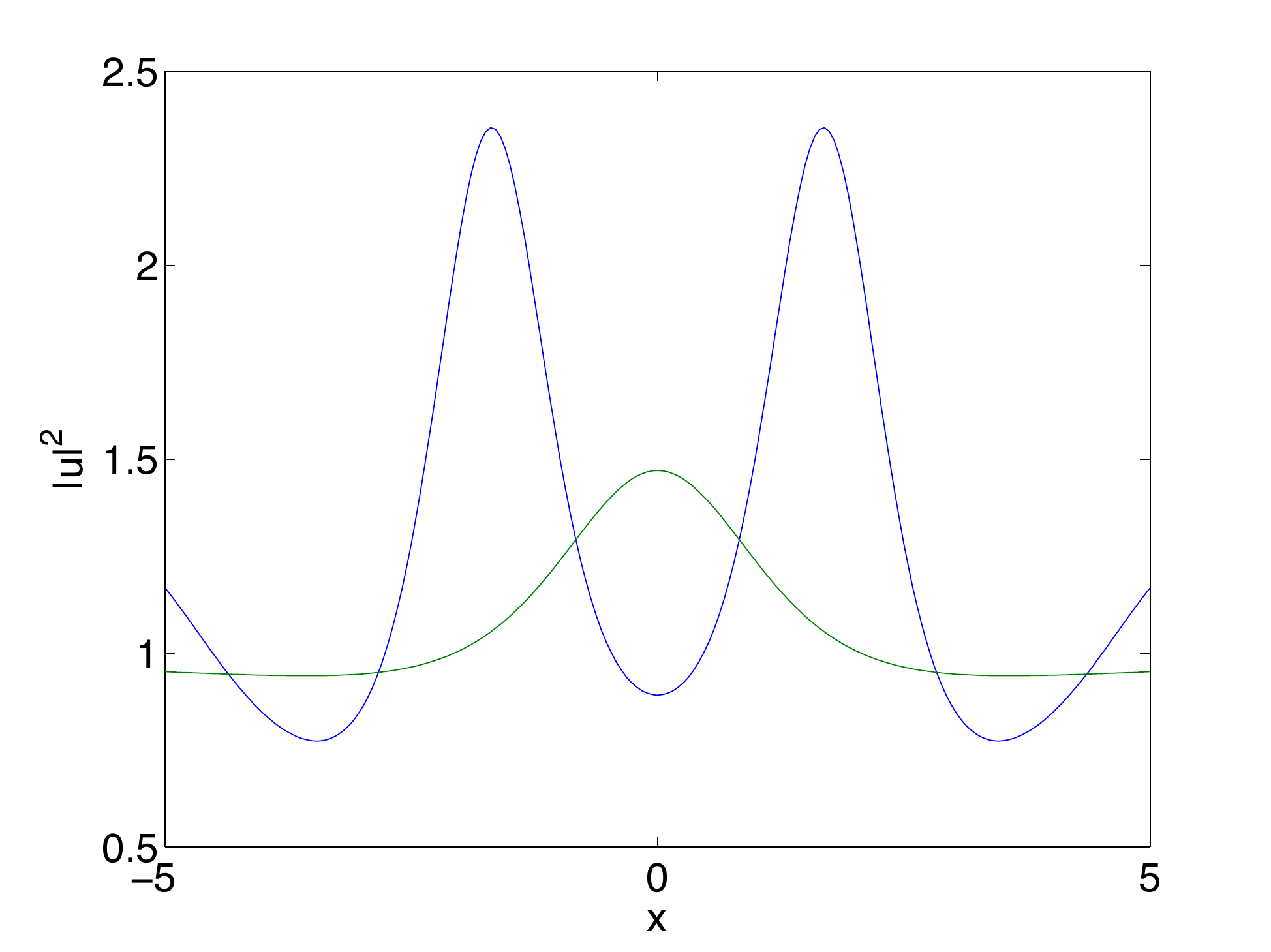}
    \includegraphics[width=0.49\textwidth]{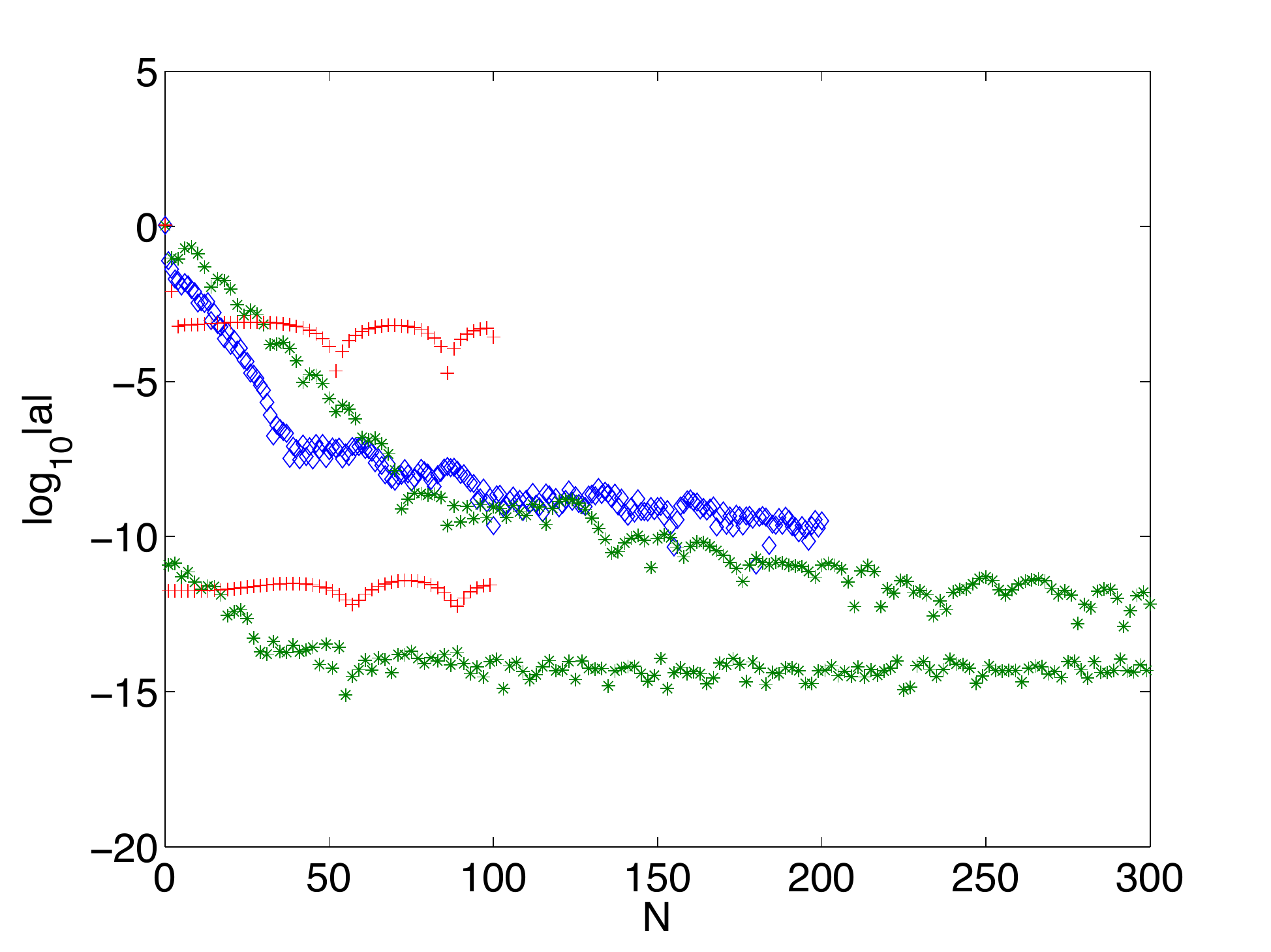}
 \caption{Solution to the NLS equation shown in 
 Fig.~\ref{breather11tm1} at $t=1$ in blue and the Peregrine breather 
 (\ref{peregrine}) for $t=1$ in green on the left, and the 
 corresponding Chebyshev coefficients for the numerical 
 solution on the right (in blue for domain I, green for domain II and 
 red for domain III).}
 \label{breather11tm1cheb}
\end{figure}

In Fig.~\ref{breather11tm1cheb} it is, however, obvious that the 
solution at the final recorded time is not close to the Peregrine 
breather. The Chebyshev coefficients in the same figure decrease to 
$10^{-5}$. During the whole computation we have 
$\Delta_{E}<2.5*10^{-3}$.

\begin{figure}[htb!]
   \includegraphics[width=0.7\textwidth]{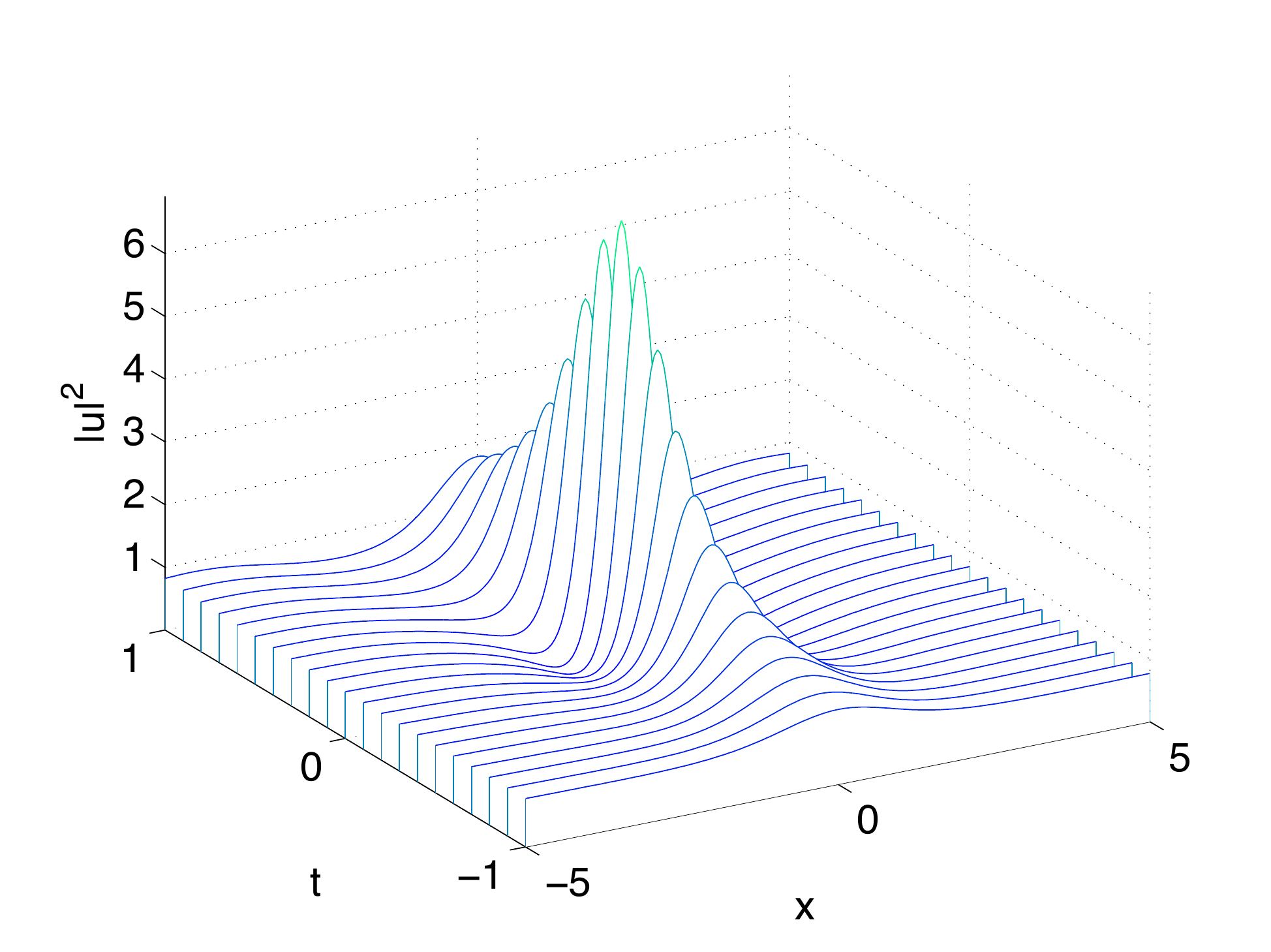}
 \caption{Solution to the NLS equation for the initial data 
 $u(x,-1)=0.9u_{Per}(x,-1)$ in dependence of time.}
 \label{breather09tm1}
\end{figure}

For the initial data $u(x,-1)=0.9u_{Per}(x,-1)$, i.e., a perturbation 
of negative $\mathcal{E}$ in (\ref{E}), the solution is shown in 
Fig.~\ref{breather09tm1}. It can be seen that the solution does not 
reach the maximum of the Peregrine breather. 

\begin{figure}[htb!]
    \includegraphics[width=0.49\textwidth]{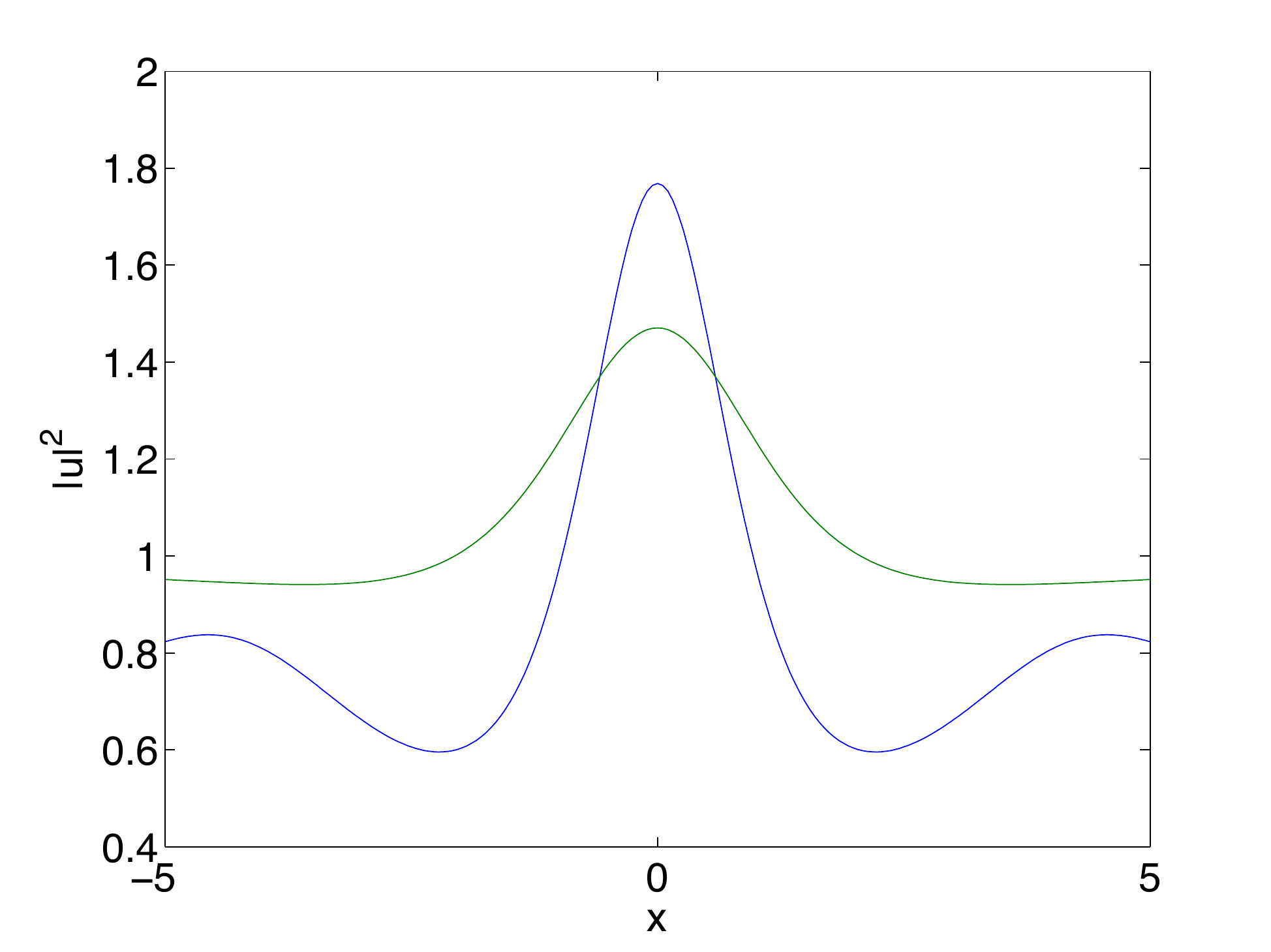}
    \includegraphics[width=0.49\textwidth]{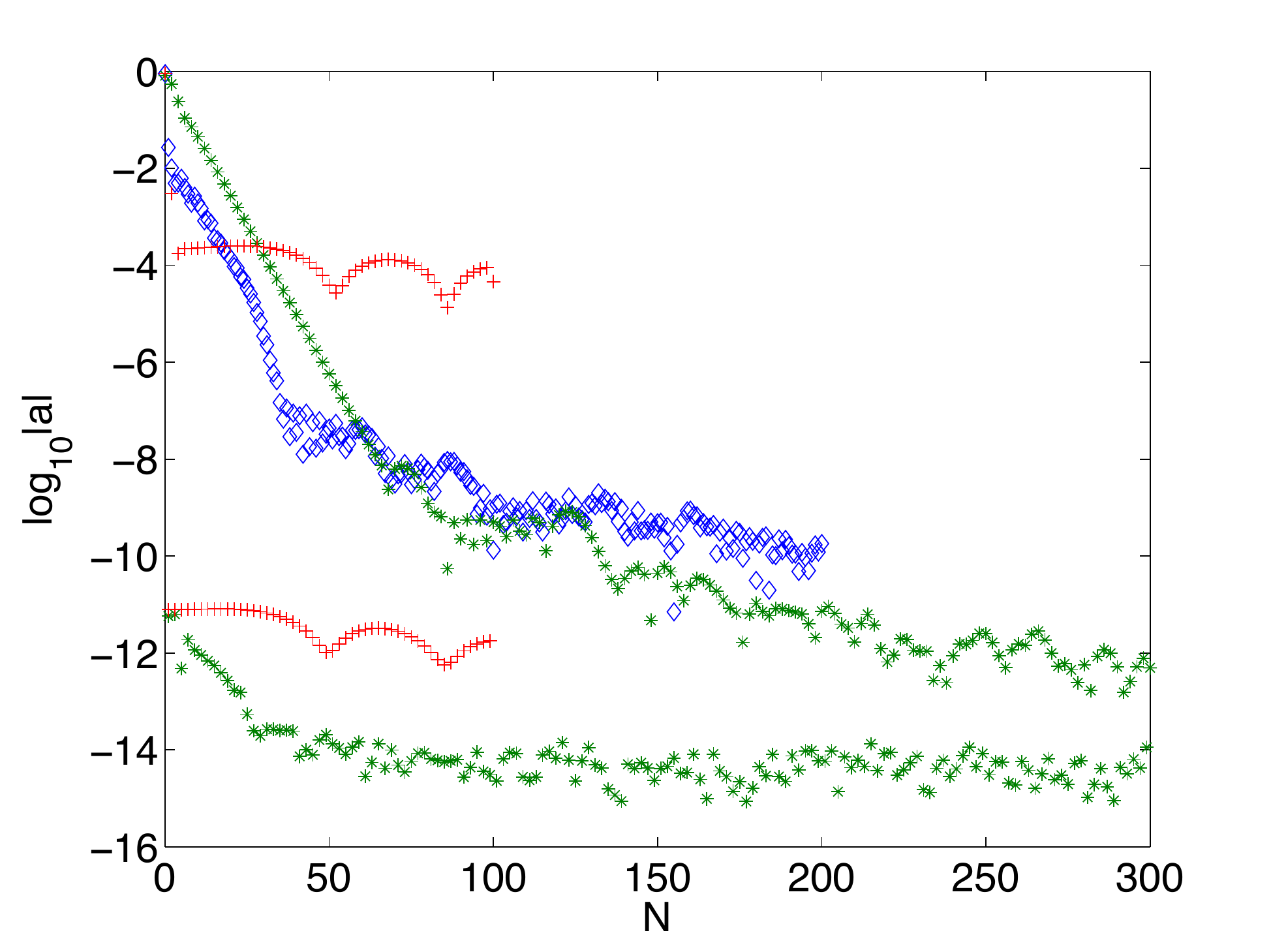}
 \caption{Solution to the NLS equation shown in 
 Fig.~\ref{breather09tm1} at $t=1$ in blue and the Peregrine breather 
 (\ref{peregrine}) for $t=1$ in green on the left, and the 
 corresponding Chebyshev coefficients for the numerical 
 solution on the right (in blue for domain I, green for domain II and 
 red for domain III).}
 \label{breather09tm1cheb}
\end{figure}

Again the solution is not close to the Peregrine breather at the 
final recorded time as can be seen in Fig.~\ref{breather09tm1cheb}. 
The solution is well resolved in the space of Chebyshev coefficients 
as can be seen in the same figure. During the whole computation we 
have $\Delta_{E}<9.05*10^{-4}$.

\section{Conclusion}
The results of the previous section indicate that solutions to the 
focusing NLS equation are not stable if they do not vanish for 
$|x|\to\infty$, and that this instability is absolute. 
This appears to be the case whether the perturbation is 
localized or not. On the other hand the focusing effect of the NLS 
equation is of course always present. In fact it seems that the 
Peregrine breather is not stable, but appears in an asymptotic sense 
explained below locally in strongly focused peaks. Note that the behavior 
shown above is actually observed in experiments in nonlinear optics 
for initial data close to the Peregrine breather \cite{kiblerp}.

\begin{figure}[htb!]
    \includegraphics[width=0.7\textwidth]{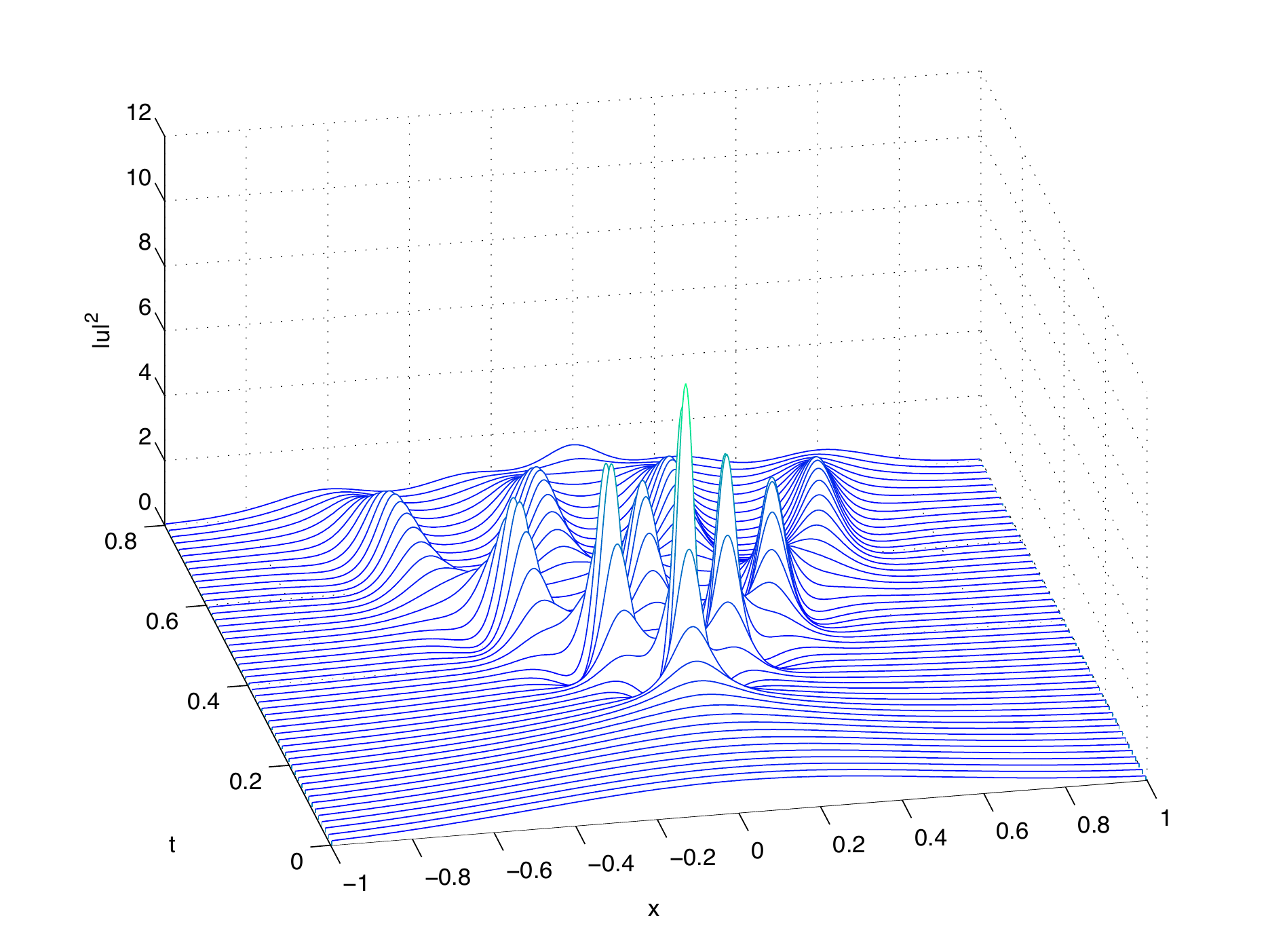}
 \caption{Solution to the NLS equation (\ref{NLSe}) with 
 $\epsilon=0.1$ for the initial data $u_{0}=\exp(-x^{2})$.}
 \label{semi}
\end{figure}

If one considers initial 
data with a support on scales of order $1/\epsilon$ and studies the 
corresponding NLS solution on time scales of order $1/\epsilon$, this 
can be conveniently done by the coordinate change $x\to\epsilon x$, 
$t\to \epsilon t$ which leads for (\ref{NLS}) to 
\begin{equation}
    i\epsilon u_{t}+\epsilon^{2}u_{xx}+2|u|^{2}u=0
    \label{NLSe}.
\end{equation}
It is well known that solutions to the NLS equation in the 
semiclassical limit $\epsilon\ll 1$ for analytic initial data in 
$L^{2}$ with a single maximum become strongly peaked at $t\sim t_{c}$ 
where $t_{c}$ is the critical time where the solution to the semiclassical system 
(the system following from (\ref{NLSe}) by formally letting 
$\epsilon\to0$) for the same initial data develops a cusp. 
An example for Gaussian 
initial data is shown in Fig.~\ref{semi}.

In \cite{DGK,DGK2} it was conjectured that the local (in the vicinity 
of the cusp) behavior of the NLS 
solution for $t\sim t_{c}$ is asymptotically for $\epsilon\to0$ given 
by the \emph{tritronqu\'ee} solution to the Painlev\'e I equation. A 
partial proof of this conjecture was given in \cite{BT}. In \cite{BT} 
it was also shown that the Peregrine breather appears in the 
asymptotic (local) description of the peak of the NLS solution in the 
semiclassical limit and is in this sense universal. 

\section*{Acknowledgement}
We thank B.~Kibler for helpful discussions.

\end{document}